%% file: compact.tex
\documentclass[11pt,amstex]{article}
\usepackage[totalwidth=480pt, totalheight=680pt]{geometry}
\usepackage[latin1]{inputenc}
\usepackage{amsmath,amssymb,amsthm}
\usepackage[dvips]{graphicx,epsfig}
\newcommand{\A}{\mathbb{A}}
\renewcommand{\AA}{\mathbb{A}}

\newcommand{\GG}{\mathbb{G}}

\newcommand{\NN}[0]{\mathbb{N}}
\newcommand{\PP}{I\!\!P}

\newcommand{\CC}{\raisebox{.6ex}{${\scriptscriptstyle /}$}\hspace{-.43em}C}
\newcommand{\ZZ}{\mbox{\rm \lower0.3pt\hbox{$\angle\!\!\!$}Z}}

\newtheorem{nt}{Notation}
\newtheorem{conjecture}[nt]{Conjecture}
\newtheorem{prop}[nt]{Proposition}

\newcommand{\fb}{\ensuremath{\downarrow}}
\newcommand{\fbd}{\ensuremath{\searrow}}

\newcommand{\fhd}{\ensuremath{\nearrow}}

\newcommand{\fd}{\ensuremath{\rightarrow}}

\newcommand{\FD}{\ensuremath{\Rightarrow}}

\newcommand{\findem}{\hfill\rule{2mm}{2mm}}

\newcommand{\fxr}{\ensuremath{\mathcal R}}

\renewcommand{\phi}{\ensuremath{\varphi}}

\newcommand{\hdx}{\ensuremath{Hilb^2(X)}}

\newcommand{\hdhdx}{\ensuremath{Hilb^2(Hilb^2(X))}}

\newcommand{\htrx}{\ensuremath{Hilb^3(X)}}

\newcommand{\htrhdx}{\ensuremath{Hilb^3(Hilb^2(X))}}

\newcommand{\inj}{\ensuremath{\hookrightarrow}}
\newcommand{\inc}{\ensuremath{\subset}}
\newcommand{\nl}{\ \\[2mm]}

\newcommand{\ox}{\otimes }

\newcommand{\plp}{\PP^2}
\renewcommand{\phi}{\ensuremath{\varphi}}

\newcommand{\s}{Spec\;}

\newcommand{\x}{\ensuremath{\times}}

\newcommand{\barr}{\begin{array}}
\newcommand{\earr}{\end{array}}
\newcommand{\bit}{\begin{itemize}}
\newcommand{\eit}{\end{itemize}}
\newcommand{\beq}{\begin{eqnarray*}}
\newcommand{\eeq}{\end{eqnarray*}}
\newcommand{\beqn}{\begin{eqnarray}}
\newcommand{\eeqn}{\end{eqnarray}}
\newcommand{\bconj}{\begin{conjecture}}
\newcommand{\econj}{\end{conjecture}}
\newcommand{\bcor}{\begin{coro}}
\newcommand{\ecor}{\end{coro} \noindent}
\newcommand{\ben}{\begin{enumerate}}
\newcommand{\een}{\end{enumerate} \noindent}
\newcommand{\bnot}{\begin{nt} }
\newcommand{\enot}{\end{nt} \noindent}
\newcommand{\bdefi}{\begin{defi}}
\newcommand{\edefi}{\end{defi} \noindent}
\newcommand{\bprop}{\begin{prop}} 
\newcommand{\brap}{\begin{rappel}}
\newcommand{\erap}{\end{rappel} \noindent }
\newcommand{\brq}{\begin{rem}}
\newcommand{\erq}{\end{rem} \noindent }
\newcommand{\bthm}{\begin{thm}}
\newcommand{\blm}{\begin{lm}}
\newcommand{\bex}{\begin{ex}}
\newcommand{\eex}{\end{ex}\noindent }
\newcommand{\bexo}{\begin{exo} \normalfont}
\newcommand{\eexo}{\end{exo}\noindent }

\font \tengothic=eufm10 scaled\magstep 1
\font\sevengothic=eufm7 scaled\magstep 1
\newfam\gothicfam
      \textfont\gothicfam=\tengothic
      \scriptfont\gothicfam=\sevengothic
\def\goth#1{{\fam\gothicfam #1}}

\newtheorem{coro}[nt]{Corollaire} 
\newtheorem{defi}[nt]{D\'efinition}

\newtheorem{ex}[nt]{Exemple} 
\newtheorem{exo}[nt]{Exercice} 
\newtheorem{lm}[nt]{Lemme} 
\newtheorem{rappel}[nt]{Rappel} 
\newtheorem{rem}[nt]{Remarque} 
\newtheorem{thm}[nt]{Th\'eor\`eme}

\newcommand{\demo}{\noindent \textit{D\'emonstration}}

\newcommand{\eprop}{\end{prop} \noindent \textit{D\'emonstration: }} 
\newcommand{\ethm}{\end{thm}\noindent \textit{D\'emonstration: }} 
\newcommand{\pdee}{{\cal P}(E)}
\newcommand{\retaprime}{R_{\eta'}(X)}
\newcommand{\elm}{\end{lm}}

\newcommand{\re}{R_{\eta}(X)}

\newcommand{\etamax}{\eta_{max}}

\newcommand{\hachun}{H_{1}(X)}
\newcommand{\hachdeux}{H_{2}(X)}
\newcommand{\hachundeux}{H_{12}(X)}
\newcommand{\hachuntrois}{H_{13}(X)}
\newcommand{\hachdeuxtrois}{H_{23}(X)}
\newcommand{\hachtrois}{H_{3}(X)}
\newcommand{\hachudt}{H_{123}(X)}

\newcommand{\feta}{F_{\eta}}

\begin{document}
\sloppy 
\title{Compactifications des espaces de configuration dans les
sch\'emas de Hilbert}
\author{L. Evain (evain@tonton.univ-angers.fr)}
\date{} 
\maketitle

\noindent
{\bf R\'esum\'e: } Soient $F(X,n):= X^n-\Delta$ 
le compl\'ementaire 
de l'union $\Delta$ des diagonales dans $X^n$,
et $U$  un quotient (\'eventuellement
trivial) de $F(X,n)$ par un sous-groupe du groupe sym\'etrique $S_n$. 
Ce travail pr\'esente des proc\'ed\'es de compactification
de $U$ dans des produits de sch\'emas de Hilbert.
Notre d\'emarche g\'en\'eralise et unifie des constructions 
classiques dues \`a  Schubert-Semple, Le Barz-Keel,
Kleiman et Cheah. 
Une \'etude g\'eom\'etrique 
plus d\'etaill\'ee est faite pour les cas 
$n\leq 3$. Cette \'etude inclut notamment une classification 
compl\`ete et une description des morphismes quotients par les actions 
naturelles. 

\section{Introduction}
 
\subsection{Quelques compactifications classiques 
et leurs applications.}
\label{les_classiques}
L'histoire des compactifications des espaces de configuration $F(X,n)$ 
est ancienne. Elle trouve son origine au si\`ecle pass\'e dans des
probl\`emes de g\'eom\'etrie \'enum\'erative. Au fil des ann\'ees, 
et jusque tr\`es r\'ecemment, de nouvelles constructions sont apparues
au gr\'e des besoins. Illustrons ces nombreuses
constructions en pr\'esentant les plus classiques. 
\nl
D\`es 1880, Schubert [Sch] utilise une 
compactification de la vari\'et\'e des triangles ($X=\plp$ et $n=3$)
pour r\'esoudre des probl\`emes 
\'enum\'eratifs. Son travail est modernis\'e par Semple
[S]. Tyrell [Ty],
Roberts et Speiser [RS], s'appuyant sur le travail de Semple,
d\'emontrent rigoureusement certaines formules de Schubert. 
Collino et Fulton [CF] 
calculent compl\`etement l'anneau d'intersection de cette compactification 
et retrouvent \'egalement les r\'esultats de Schubert. 
Le Barz \'etend la construction de 
Schubert-Semple \`a toute vari\'et\'e lisse par l'utilisation de sch\'emas
de Hilbert [LB] et Keel l'\'etend \`a tout sch\'ema par une approche 
fonctorielle [Kee]. \\
Kleiman [K] construit une compactification de $F(X,n)$
par r\'ecurrence sur $n$ en utilisant des points infiniment voisins.
Cette approche lui permet d'obtenir des formules d\'ecrivant le lieu multiple
d'un morphisme. Dolgachev et Ortland [DO] en d\'eduisent d'autres
constructions en liaison avec les fonctions theta. 
\\
La compactification de Fulton-MacPherson [FMP] admet plusieurs 
d\'efinitions, soit fonctorielle, soit 
g\'eom\'etrique, l'une d'entre elles \'etant un \'eclatement 
subtil de $X^n$. Quand $X$ est compact, elle  permet le calcul du type d'homotopie 
rationnel des espaces de configuration 
en fonction des invariants de $X$.
\\
Citons enfin le travail de Cheah [Ch] qui a \'etudi\'e des 
compactifications 
similaires dans l'esprit 
\`a celle de LeBarz,  se  projetant sur le sch\'ema de Hilbert $Hilb^n(X)$.

\subsection{Le probl\`eme}
Expliquons plus en d\'etail l'approche de Le Barz, qui est le point de
d\'epart de notre travail. Etant donn\'es trois points distincts $p_1,p_2,p_3$
de $X$, on peut former les trois doublets $d_{12}:=p_1\cup p_2$, $d_{13},
d_{23} \in Hilb^2(X)$ et le triplet $p_{123}=p_1 \cup p_2 \cup p_3
 \in Hilb^3(X)$. Cette construction se reformule en disant que 
$F(X,3)$ est isomorphe \`a un sous-sch\'ema $Z(X)$ localement ferm\'e de 
$X^3 \x Hilb^2(X)^3 \x Hilb^3(X)$. L'adh\'erence $\overline {Z(X)}$ est 
donc une compactification naturelle de $F(X,3)$ dans un produit de 
sch\'emas de Hilbert. Le Barz a montr\'e que cette adh\'erence, 
\`a priori difficilement manipulable, peut en fait \^etre d\'ecrite 
g\'eom\'etriquement en termes de lieux d'incidence: les points de 
l'adh\'erence sont les $7$-uplets $(p_1,p_2,\dots,p_{123})$ satisfaisant
les relations \'evidentes 
$p_1\inc d_{12}\inc p_{123}$, $p_3$ est le r\'esiduel de $d_{12}$ dans 
$p_{123}$ et les relations s'en d\'eduisant par sym\'etrie.
Keel a remarqu\'e que cette description par adh\'erence pouvait 
\^etre exploit\'ee pour donner une d\'efinition de $\overline {Z(X)}$ 
comme repr\'esentant d'un certain foncteur et en d\'eduire quelques
cons\'equences g\'eom\'etriques. On obtient donc 
finalement une compactification agr\'eable \`a manipuler car elle
jouit d'une triple d\'efinition, par adh\'erence, par incidence, 
et fonctorielle. Cette multiplicit\'e des points de vue est tr\`es 
semblable \`a l'approche de [FMP]. 
\\
Dans la fronti\`ere de la compactification construite par Le Barz, 
toutes les informations sur la collision des points $p_1,p_2,p_3$ se 
trouvant dans les sch\'emas de Hilbert n'ont visiblement pas \'et\'e 
exploit\'ees. Par exemple, deux doublets
$d_{12}$ et $d_{13}$ de $Hilb^2(X)$ distincts
d\'efinissent un point $d^1$ de $Hilb^2(Hilb^2(X))$. Utilisant
les points $d^1,d^2,d^3$, on peut construire par adh\'erence une 
compactification de $F(X,3)$ dans un produit plus gros contenant 
des facteurs de la forme $Hilb^2(Hilb^2(X))$. R\'eit\'erant
le processus, on peut construire des compactifications dans des 
espaces produits $P_1 \x \dots \x P_n$, chacun des termes $P_i$ 
du produit \'etant un sch\'ema de Hilbert embo\^\i t\'e 
$Hilb^{p_l}(Hilb^{p_{l-1}}(\dots Hilb^{p_1}(X)))$. En variant  les termes 
$P_i$, on a donc une infinit\'e de compactifications \`a notre 
disposition. 
\nl
Le travail qu'on se propose d'effectuer ici est 
d'\'etudier et de classifier 
les compactifications obtenues par ce proc\'ed\'e.

\subsection{Les r\'esultats}
Cet article s'articule en trois parties. 
Premi\`erement, on
d\'egage les notions d'incidence et les foncteurs qui permettent 
de manipuler les compactifications ais\'ement via une triple 
d\'efinition (par adh\'erence, par incidence et fonctorielle - section
\ref{sec:definitions}). 
Ensuite, on classifie les compactifications de $F(X,n)$ et de 
ses quotients obtenus quand $n\leq 3$ et on \'etudie la g\'eom\'etrie 
des compactifications obtenues (sections \ref{sec:classification}
et \ref{sec:etude_des_quotients}). Enfin, on compare les compactifications
obtenues aux compactifications classiques d\'ecrites ci-dessus 
(section \ref{sec identification des classiques}).
\\
Le cas $n=2$ \'etant trivial, on ne pr\'esente dans cette introduction
que les r\'esultats concernant les  compactifications de $F(X,n=3)$. 
On peut extraire de cette \'etude deux r\'esultats surprenants.
Tout d'abord, alors qu'il existe une infinit\'e de choix possibles 
pour le produit dans lequel on construit la compactification 
de $F(X,3)$, il n'existe \`a post\'eriori qu'un nombre fini de
compactifications \`a isomorphisme pr\`es. 
\begin{thm} \label{thm: classification}
Il y a \`a isomorphisme de compactification pr\`es 11
compactifications de $F(X,3)$ o\`u de ses quotients obtenues par 
notre proc\'ed\'e.De plus, toute compactification
est isomorphe \`a une compactification dans $P_1\x \dots \x P_r$ 
o\`u chaque $P_i$ est soit de la forme $Hilb^i(X)$, soit de la 
forme $Hilb^i(Hilb^j(X))$. 
\end {thm}
\noindent
Ces vari\'et\'es seront not\'ees 
$R_{123}(X),R_{123}^1(X)...$  La signification pr\'ecise des notations 
sera d\'evelopp\'ee dans le corps du texte et
n'est pas n\'ecessaire \`a  cette introduction
( grossi\`erement, les indices sont relatifs aux termes 
$Hilb^i(X)$ tandis que les exposants sont relatifs aux termes 
$Hilb^i(Hilb^j(X))$. Par exemple, la vari\'et\'e $R^1_{123}$ est
l'adh\'erence du morphisme $f:F(X,3) \fd Hilb^3(X)\x Hilb^2(Hilb^2(X))$ 
d\'efini par $f(p_1,p_2,p_3)=(p_{123},d^1)$ avec les notations 
pr\'ec\'edentes).
\nl
Le deuxi\`eme fait inattendu est que les structures d'ordre sup\'erieur
(c'est \`a dire 
les termes de la forme 
$Hilb^{p_l}(Hilb^{p_{l-1}}(\dots
Hilb^{p_1}(X)))$,
avec $l\geq 2$)
peuvent \^etre utilis\'ees pour d\'ecrire les passages au
quotient. Par exemple, la vari\'et\'e construite par Le Barz
dans $X^3 \x Hilb^2(X)^3 \x Hilb^3(X)$ est 
naturellement munie d'une action du groupe sym\'etrique 
$S_3$ et la structure d'ordre deux
$Hilb^3(Hilb^2(X))$
permet une description explicite du quotient: 

\begin{thm} 
Le quotient de la vari\'et\'e de Le Barz est isomorphe \`a 
l'adh\'erence de $F(X,3)/S_3$ dans $Hilb^3(Hilb^2(X))\x Hilb^3(X)$. 
\end {thm}
\noindent
Ce th\'eor\`eme \'etait en fait notre motivation premi\`ere: obtenir
une description du quotient de la vari\'et\'e de Le Barz ais\'ement
manipulable via ses propri\'et\'es universelles (voir 
le contexte \`a la fin de cette introduction). 
En d'autres termes, les structures de niveau sup\'erieur
s'imposent d'elles m\^emes quand il s'agit d'\'etudier les quotients. 
Ce principe est illustr\'e par le th\'eor\`eme g\'en\'eral suivant
qui dit que les compactifications construites 
ont le m\'erite de former une classe 
de compactifications stable par quotient.

\begin{thm} \label{thm: quotients}
Supposons le corps $k$ de caract\'eristique diff\'erente de deux et trois.
Les groupes agissant sur 
$$R_{123}(X),\ R_{1,123}(X),\ R_{1,2,3,12,123}(X),\  R_{3,12,123}(X)$$ 
$$ 
R^1_{123}(X),\ R^1_{1,2,3,12,13,123}(X),\ R^{123}_{123}(X),\ R^{123}_{1,123}(X)$$  
$$ 
\ R^{3,123}_{3,12,123}(X),\ R^{3,123}_{1,2,3,12,123}(X) ,
\  R_{max}(X)$$
sont 
les groupes sym\'etriques
$$
S_1,\ S_1,\ S_2,\ S_1$$
$$
  \ S_1,\ S_2,\ S_1,\ S_1$$  
$$
 \ S_1,\ S_2\ S_3 $$
et les quotients respectifs sont 
$$ 
R_{123}(X),\  R_{1,123}(X),\  R_{3,12,123}(X),\ R_{3,12,123}(X)$$ 
$$ 
 R^{1}_{123},\  R^{1}_{123},\  R^{123}_{123}(X),\ R^{123}_{1,123}(X)$$ 
$$ 
R^{3,123}_{3,12,123}(X),\  R^{3,123}_{3,12,123}(X),\  R^{123}_{123}(X)
$$
Le quotient partiel de 
$ R_{max}(X) $
par $S_2$ est  
$ R^{3,123}_{3,12,123}(X)$.
\end{thm}

\noindent
De fa\c con \'evidente, si $P_1$ et $P_2$ sont des produits 
de sch\'ema de Hilbert embo\^\i t\'es, l'adh\'erence de 
$F(X,3)$ dans $P_1\x P_2$ se projette sur l'adh\'erence dans 
$P_1$. En termes fonctoriels, on a un morphisme d'oubli.
Le th\'eor\`eme suivant dit que les morphismes quotients 
du dernier \'enonc\'e sont essentiellement des morphismes d'oubli.
 
\begin{thm} \label{thm: quotient=oubli}
  Soient $R_{\alpha}$ et $R_{\beta}$ deux compactifications de $F(X,3)$
telles que $R_{\beta}$ soit le quotient de $R_{\alpha}$ par un sous-groupe 
$G$ du groupe sym\'etrique $S_3$. Alors il existe une 
compactification $R_{\alpha}'$ isomorphe \`a $R_{\alpha}$ telle que le 
passage au quotient $R_{\alpha}' \fd R_{\beta}$ soit un morphisme d'oubli.
\end{thm}

\noindent
Les rapports entre les diff\'erentes compactifications
est r\'esum\'e par le th\'eor\`eme suivant, \'etabli
au cours de la d\'emonstration de la classification, et 
cl\'e des autres th\'eor\`emes. 

\bthm \label{thm: stratification et oubli compatibles}
Chacune des compactifications est munie d'une stratification 
naturelle. Pour chacune des stratifications, il existe une 
strate g\'en\'erale, ouverte et dense, et une strate sp\'eciale 
incluse dans l'adh\'erence de toutes les autres strates. Les 
morphismes d'oubli liant les compactifications sont donn\'es par 
les fl\`eches du diagramme suivant et leurs compositions. 

  \begin{figure}[h] 
     \begin{center}
        \input{mor_oublis.pstex_t}
     \end{center} 
  \end{figure}
De plus, l'image inverse d'une strate par l'un quelconque de 
ces morphismes est 
une r\'eunion de strates. 
\end{thm}

\noindent
Enfin,
les liens unissant les vari\'et\'es 
que nous avons construites et les constructions classiques 
sont explor\'es et r\'esum\'es par le th\'eor\`eme suivant. 
\begin{thm} \label{comparaison aux classiques}
  Les vari\'et\'es
  $R_{max}(X),\ R_{1,2,3,12,13,123}^1(X),\ R_{1,2,3,12,123}(X),$
  $\ R_{1,23,123}^{1,123}(X),$ 
  $\ R_{123}^{123}(X)$
  $\ R^1_{123}(X)$
  sont respectivement les  vari\'et\'es de LeBarz $\widehat{H_3}(X)$, 
  de Kleiman $K_3(X)$, 
  de Cheah, les quotients  $\widehat{H_3}(X)/S_2$,
  $\widehat{H_3}(X)/S_3$, $K_3(X)/S_2$. En outre, $R_{max}(X)$
s'identifie \`a la vari\'et\'e de Schubert Semple dans le cas o\`u 
$X=\plp$. 
\end{thm}

\noindent
Il convient de 
pr\'eciser la signification du terme ``compactification''
figurant dans les \'enonc\'es pr\'e\-c\'e\-dents. 
Comme expliqu\'e,
les vari\'et\'es  consid\'er\'ees  vivent 
dans des produits de sch\'emas de Hilbert embo\^\i t\'es
$Hilb^{p_l}(Hilb^{p_{l-1}}(((...(Hilb^{p_1}(X)))))$ et  
admettent trois d\'efinitions diff\'erentes: comme adh\'e\-rence de 
morphisme, comme lieu sch\'ematique d'incidence, ou 
par repr\'esentabilit\'e de certains foncteurs. 
Si on se fixe 
un produit $P$ de sch\'emas de Hilbert embo\^\i t\'es, on disposera 
d'un morphisme $F(X,n) \fd P$ et les compactifications qui nous 
int\'eressent sont les adh\'erences des images.
Dans le cas o\`u $P=Hilb^3(Hilb^2(X))\x Hilb^3(X)$, 
on peut donner une d\'efinition satisfaisante du lieu sch\'ematique 
 param\'etrant les couples $(d,t)$ 
pour lesquels le triplet de doublets $d$ est inclus dans le triplet $t$.
On montre que ce lieu d'incidence  co\"\i ncide avec 
l'adh\'e\-rence du morphisme $F(X,3)\fd P$.  
La situation g\'en\'erale est un peu moins 
simple: \'etant donn\'es un produit $P$
de sch\'emas de Hilbert embo\^\i t\'es et le morphisme $F(X,3)\fd P$
correspondant, il y aura un lieu sch\'ematique 
d'incidence naturel $R\subset P$ tel que l'adh\'erence $A$ de
l'image v\'erifie $A \inc R$. 
Le sch\'ema $R$ repr\'esente toujours un foncteur 
mais on n'aura pas en g\'eneral $A=R$. 
Les ``belles'' compactifications, celles que nous chercherons 
\`a classifier, seront celles pour lesquelles $A=R$. C'est le sens du
mot compactification dans les \'enonc\'es.

\subsection{Application de ces constructions}
A l'instar des compactifications classiques pr\'esent\'ees au debut de
cette introduction, il existe des applications
ayant motiv\'e nos constructions.
M\^eme si ces applications
seront d\'evolopp\'ees ailleurs, nous les donnons ici
\`a titre de motivation. 
\\
Le probl\`eme original consistait 
\`a comprendre les collisions de trois gros points sur une surface 
lisse $S$, o\`u rappelons le, un gros point de taille $m$
est un sous-sch\'ema d\'efini par la puissance $m^{eme}$ d'un 
id\'eal maximal. Plus pr\'ecis\'ement, notons $Coll(n_1,n_2,n_3)(S)$
la sous-vari\'et\'e irr\'eductible de $Hilb(S)$ dont le point 
g\'en\'erique param\`etre la r\'eunion g\'en\'erique de trois gros 
points de $S$ de taille $n_1,n_2,n_3$. Les collisions de trois 
gros points sont les points de la fronti\`ere de 
$Coll(n_1,n_2,n_3)(S)$. Le r\'esultat principal de [Ev] dit 
que lorsque $n_1,n_2,n_3$ parcourent $\NN^3$, alors 
$Coll(n_1,n_2,n_3)(S)$ ne parcourt qu'un nombre fini de classes 
d'isomorphismes. 
La m\'ethode consiste \`a construire des isomorphismes explicites 
avec certaines des compactifications d\'ecrites  dans le 
pr\'esent article. On en d\'eduit une classification des collisions 
en \'etudiant la restriction des isomorphismes
\`a leur  fronti\`ere.
\nl
Citons pour finir deux probl\`emes classiques pour lesquels les 
compactifications par sch\'emas  de Hilbert embo\^\i t\'es pourraient se
r\'ev\'eler utiles. 
Si $X$ est une vari\'et\'e lisse, les probl\`emes de construire de
fa\c con explicite
une compactification lisse de $F(X,n)$ se projetant sur le sch\'ema de
Hilbert $Hilb^n(X)$ et une compactification lisse de $X^n/S_n$ 
restent ouverts.
On peut esp\'erer 
une r\'eponse positive \`a la question suivante:
existe-t-il des compactifications construites dans des produits de
sch\'emas de Hilbert embo\^\i t\'es qui soient solutions 
des deux probl\`emes pr\'ec\'edents ?
\\
Il r\'esulte de ce travail que la r\'eponse est oui pour $n=3$
(les compactifications qui conviennent sont 
$R_{max}(X)$ et $R_{123}^{123}(X)$).

\nl
Je remercie vivement A. Hirschowitz pour ses conseils 
lors de la r\'ealisation 
de ce travail.

\section{D\'efinition des compactifications}
\label{sec:definitions}
Dor\'enavant, $X$ est un sch\'ema projectif sur un corps 
alg\'ebriquement clos $k$ de caract\'eristique quelconque
(quasi-projectif conviendrait \'egalement en adaptant 
quelques d\'emonstrations). 
\nl
Dans cette section, on d\'efinit les compactifications et on 
donne les premi\`eres propri\'et\'es d\'ecoulant directement 
des d\'efinitions. Dans la section \ref{subsec:def par adherence}, on 
d\'efinit un morphisme $f_{\eta}:F(X,n)\fd H_{\eta}$ o\`u $H_{\eta}$ est 
un produit de sch\'emas de Hilbert embo\^\i t\'es d\'ependant 
d'une donn\'ee combinatoire $\eta$. Pour chaque $\eta$, il existe 
un sous-groupe $G_{\eta}$ du groupe sym\'etrique $S_n$
tel que l'adh\'erence $A_{\eta}:=
\overline{f_{\eta}(F(X,n))}$
soit une compactification du quotient $F(X,n)/G_{\eta}$.
Dans la section \ref{sec:def par incidence}, 
on d\'efinit \`a l'aide de relations 
d'incidence un lieu $R_{\eta}$ dans $H_{\eta}$ tel que $A_{\eta}\inc 
R_{\eta}$. Dans la section \ref{subsec:def par foncteurs}, 
on introduit des foncteurs 
$F_{\eta}$ dont on montre qu'ils sont repr\'esentables par 
$R_{\eta}$. 
\\
Les compactifications $A_{\eta}$ qui nous int\'eressent sont celles 
pour lesquelles $A_{\eta}=R_{\eta}$. Autrement dit, ce sont celles 
dont la structure est riche de sorte qu'elles puissent \^etre 
d\'efinies au choix par adh\'erence (via $f_{\eta}$), par
lieu sch\'ematique d'incidence (via $R_{\eta}$) ou encore 
par fonctorialit\'e (via $F_{\eta}$).

\subsection{D\'efinition par adh\'erence} 
\label{subsec:def par adherence}
On va d\'efinir successivement un enrichissement $\eta$,
le sch\'ema de Hilbert associ\'e $H_{\eta}$, le morphisme 
$f_{\eta}:F(X,n)\fd H_{\eta}$, le groupe $G_{\eta}$. 
La section se conclut par la proposition 
\ref{description de l'adherence comme quotient}
qui dit que $A_{\eta}:=\overline{f_{\eta}(F(X,n))}$ est une compactification 
du quotient $F(X,n)/G_{\eta}$. 
\nl
Soit $E=\{1,\dots,n\}$. Pour $p>0$, notons $\Sigma_p(E)$ 
les sous-ensembles de cardinal $p$
de $E$.
D\'efinissons r\'ecursivement 
$\Sigma_{p_l,...p_1}(E):=\Sigma_{p_l,...p_{2}}(\Sigma_{p_1}(E))$
et $\Sigma(E)=\cup_{l \in \NN^*, (p_l,\dots,p_1)\in (\NN^{*})^{l}}
\Sigma_{p_l,...p_1}(E)$. 
Un \'el\'ement de $\Sigma(E)$ (resp. de $\Sigma_{p_l,...p_1}(E)$ )
est appel\'e une structure
sur $E$ (resp. une structure 
de niveau $l$). Notons $Hilb^p(X)$ le sch\'ema de Hilbert param\'etrant
les sous-sch\'emas ponctuels de $X$ de longueur $p$.
On pose 
$Hilb^{p_l,p_{l-1},\dots,p_1}(X):=
Hilb^{p_l}(Hilb^{p_{l-1}}(\dots (Hilb^{p_1}(X))))$.
Si $\sigma \in \Sigma_{p_l,...p_1}(E)$, on posera $H_{\sigma}(X)=
Hilb^{p_l,p_{l-1},\dots,p_1}(X)$ et$H_{\sigma}^-(X)=
Hilb^{p_{l-1},\dots,p_1}(X)$ .
Pour tout \'el\'ement  $\sigma$ de $\Sigma_{p_l,...p_1}(E)$ (resp. de 
$Hilb^{p_l,p_{l-1},\dots,p_1}(X)$), on note 
$[\sigma]$ le sous-ensemble de $\Sigma_{p_{l-1},...p_1}(E)$
(resp. le sous-sch\'ema 
de $Hilb^{p_{l-1},\dots,p_1}(X)$) param\'etr\'e par $\sigma$.

\bdefi
Soit $\sigma \in \Sigma(E)$. On note 
$
f_{\sigma}:
F(X,n){\fd} H_{\sigma}(X)$ le morphisme qui envoie $x=(x_1,\dots,x_n)$ sur
$f_{\sigma}(x)$ v\'erifiant:
\bit
\item
$[f_{\sigma}(x)]=\cup_{j\in [\sigma]}\;x_j$  si $\sigma$ est de 
niveau un
\item
$[f_{\sigma}(x)]=f_{\sigma_1}(x) \cup \dots \cup f_{\sigma_q}(x)$
si $\sigma$ est de
niveau sup\'erieur \`a un avec  $[\sigma]=\{\sigma_1,\dots,\sigma_q\}$.
\eit
\edefi

\newcommand{\ssig}{\Sigma_{p_l,p_{l-1},\dots,p_1}(E)}
\brq
\label{identification}
Il existe une bijection canonique entre $\ssig$ et 
$\Sigma_{p_l,p_{l-1},\dots,p_r,1,p_{r-1},\dots,,p_1}(E)$, et 
entre $Hilb^{p_l,p_{l-1},\dots,p_1}(X)$ et 
$Hilb^{p_l,p_{l-1},\dots,p_r,1,p_{r-1},\dots,,p_1}(X)$.
Dans la suite, nous identifierons 
deux \'el\'ements qui se correspondent par l'une de ces  bijections.
\erq
On appelle enrichissement de $E$ (ou enrichissement tout 
court lorsque le contexte est clair) un ensemble 
$\eta=\{\sigma_1,\dots,\sigma_s\}$ de structures de 
$\{1,\dots,n\}$ contenant la structure $\sigma=\{1,\dots,n\}$ 
de $\Sigma_n(E)$. Pour des raisons
li\'ees \`a la description des  
actions de groupe, on n'identifiera pas 
deux enrichissements $\{\sigma_1,\dots,\sigma_s\}$
et $\{\sigma_{p(1)},\dots,\sigma_{p(s)}\}$ 
qui diff\`erent l'un de l'autre par une permutation $p$
des facteurs. Le niveau d'un enrichissement sera le plus 
grand des niveaux des $\sigma_i$.

\bdefi
Soit $\eta=\{\sigma_1,\dots,\sigma_s\}$
un  enrichissement. 
On d\'efinit
$H_{\eta}(X):=H_{\sigma_1}(X)\x H_{\sigma_2}(X)\x\dots \x H_{\sigma_s}(X)$
et 
$f_{\eta}:=f_{\sigma_1}\x f_{\sigma_2}\x\dots \x f_{\sigma_s}
:F(X,n)\fd H_{\eta}(X)$.
\edefi

Le groupe $S(E)=S_n$ agit sur 
$X^n$ par $\tau.(x_1,\dots,x_n)=(x_{\tau(1)},\dots,x_{\tau(n)})$.
L'action naturelle de $S(E)$ sur $E$ 
induit une action de $S(E)$ sur 
l'ensemble $\Sigma(E)$ des enrichissements. 
Un enrichissement $\eta$ \'etant fix\'e, 
on note $G_{\eta}$ le sous-groupe de $S(E)$ stabilisant $\eta$.

\begin{prop} \label{description de l'adherence comme quotient}
Le plus petit sous-sch\'ema ferm\'e  $\overline{f_{\eta}(F(X,n))}$
de $H_{\eta}(X)$ qui factorise  $f_{\eta}$ est une compactification 
du quotient  $F(X,n)/G_{\eta}$ qui se projette sur $Hilb^n(X)$. 
\end{prop}
\noindent
\textit{D\'emonstration}: puisque $\eta$ contient l'enrichissement 
$\{1,\dots,n\}$ par d\'efinition, $Hilb^n(X)$ est un facteur 
de $H_{\eta}(X)$ et l'existence d'une projection 
de $\overline{f_{\eta}(F(X,n))} \inc H_{\eta}(X)$ 
sur $Hilb^n(X)$
en r\'esulte. Montrons que deux \'el\'ements ${x}$ 
et ${y}$ de $F(X,n)$ ont m\^eme image par
$f_{\eta}$ ssi $\exists g \in G_{\eta}$ t.q. $g.x=y$. 
Si  $x$ et $y$ ont m\^eme image par $f_{\eta}$, 
alors ils ont m\^eme image par 
 la 
compos\'ee $F(X,n)\fd H_{\eta}(X)\fd Hilb^n(X)$, donc 
il existe $g\in S(E)$ tel que $g.x=y$.
Pour conclure, il suffit donc de montrer pour $g \in S(E)$ 
que $f_{\eta}(x)=f_{\eta}(g.x)$ ssi $g \in G_{\eta}$.
Il suffit pour cela de remarquer que  
$f_{\eta}(g.x)=f_{g.\eta}(x)$ et que,  $x$ \'etant fix\'e
dans $F(X,n)$, $f_{\eta}(x)=f_{\eta'}(x)$ ssi $\eta=\eta'$
(r\'ecurrences faciles sur 
le niveau de $\eta$).
\findem

\begin{rem}
Dans la suite, on emploiera simplement le mot ``compactification''
sans
pr\'eciser de quel quotient de $F(X,n)$ il s'agit.  
\end{rem}

\begin{rem}
On pourrait d\'efinir un enrichissement $\eta$ sans imposer la 
condition: $\eta$ contient
$\{1,\dots,n\}$. 
Il faudrait alors demander que la projection naturelle de 
$\Sigma(E)$ dans ${\cal P}(E)$ envoie $\eta$ sur $E$ pour 
que le morphisme $f_{\eta}$ soit \`a fibres finies et que l'image 
soit un quotient de $F(X,n)$ par un sous-groupe de $S(E)$. 
On obtiendrait avec cette d\'efinition plus de compactifications
mais les nouvelles compactifications
ne se projetteraient plus sur $Hilb^n(X)$. 
\end{rem}

\subsection{D\'efinition par incidence.}
\label{sec:def par incidence} 
On suppose ici fix\'e un enrichissement $\eta$ de $E$, et on a donc 
un morphisme $f_{\eta}:F(X,n)\fd H_{\eta}(X)$ dont on veut caract\'eriser 
g\'eom\'etriquement l'adh\'erence en termes d'incidence. 
Pour cela, on d\'efinit dans cette section 
les incidences liant les structures de 
$\eta$. Formellement, on pose $I(\eta')=1$ si 
$\eta'=\{\sigma,\sigma_1,\dots,\sigma_s\}\inc \eta$
v\'erifie $\sigma\inc \sigma_1 \cup  \dots \cup \sigma_s$.
Pour chaque sous-enrichissement $\eta'\inc \eta$ tel que 
$I(\eta')=1$, on d\'efinit un sous-sch\'ema $Inc_{\eta',\eta}$ de 
$H_{\eta}(X)$ et on pose $R_{\eta}(X):=\cap _{I(\eta')=1}Inc_{\eta',\eta}$. 
On v\'erifie que pour $X$ lisse irreductible, 
$f_{\eta}:F(X,n)\fd H_{\eta}(X)$
se factorise par $R_{\eta}(X)$.
Le sch\'ema $R_{\eta}(X)$ est le lieu d'incidence naturel de $H_{\eta}(X)$ 
pour lequel on esp\`ere avoir l'\'egalit\'e $H_{\eta}(X)=\overline{f_{\eta}(F(X,n))}$.
Malheureusement, on verra que m\^eme en 
mettant toutes les inclusions possibles, l'adh\'erence et le 
lieu d'incidence ne co\"\i ncident pas toujours. Un enrichissement 
pour lequel l'adh\'erence et le lieu d'incidence co\"\i ncident 
pour toute vari\'et\'e lisse irr\'eductible $X$ sera  appel\'e enrichissement 
admissible.

\subsubsection{D\'efinition des relations d'incidence sur les
  structures} 
\label{sssec: def_incidence_structure}
Soit $\eta'=\{\sigma,\sigma_1,\dots,\sigma_s\}$ un enrichissement. 
D\'efinissons le bool\'een $I(\eta')$. 
On pose $I(\eta')=0$ s'il n'existe pas d'entiers 
$p_1,\dots,p_l$, $n_1, \dots,n_s$, $q_1,\dots,q_r$ tels que 
$\sigma\in \Sigma_{q_r,\dots,q_1,p_l,\dots,p_1}$, 
et $\sigma_i \in \Sigma_{n_i,p_l,\dots,p_1}$. Si de tels entiers
existent, on d\'efinit $I(\eta')$ par r\'ecurrence sur $r$. 
Pour $r=1$, on pose $I(\eta')=1$ ssi $[\sigma]\inc [\sigma_1]\cup
\dots \cup [\sigma_s]$ dans $\Sigma_{p_l,\dots,p_1}(E)$. 
Pour $r>1$ et  $\sigma=\{\tau_1,\dots,\tau_{q_r}\}$, on 
pose $I(\eta')=1$ ssi $\forall i,\ I(\{\tau_i,
\sigma_1,\dots,\sigma_s\})=1$. 

\subsubsection{D\'efinition des relations d'incidence sur les
sch\'emas}
Soit $P$ un sch\'ema projectif, $B$ un sch\'ema quelconque,
$F$ et $G$ deux sous-sch\'emas ferm\'es respectifs de 
$Hilb^{q_r,\dots,q_1}(P)\x B$ et $P\x B$.
\\
Pour $\eta'$ tel que $I(\eta')=1$, on reprend les notations 
de \ref{sssec: def_incidence_structure} et on d\'efinit un 
sous-sch\'ema $Inc_{\eta'}$ de $H_{\eta'}(X)$ par 
r\'ecurrence sur $r$.
\\
\textit{Le cas r=1}. 
Rappelons la proposition suivante ([Kee]):
\begin{prop}
\label{def de Inc}
  Soient  $F\fd B$ et $G\fd B$ deux morphismes 
  tels que $F\fd B$ soit plat et fini. Il existe un plus 
  grand sous-sch\'ema ferm\'e $Inc_B(F,G)$ de $B$ tel que $F\x_B Inc(F,G)
  \inc G\x_B Inc(F,G)$. En outre, $Inc_B(F,G)$ peut \^etre 
caract\'eris\'e par la propri\'et\'e suivante: un morphisme 
$\phi: Z \fd B$ se factorise par $Inc_B(F,G)$ ssi $F\x _B Z \inc G\x
  _B Z$.
\end{prop}
\noindent
Soient $U,U_1, \dots, U_s$ les ferm\'es universels de
$$(Hilb^{q_1,p_l,\dots,p_1}(X)\x \Pi_{i=1}^{i=s} 
Hilb^{n_i,p_l,\dots,p_1}(X))\x Hilb^{p_l,\dots,p_1}(X)$$
dont les fibres respectives au dessus de 
$(h,h_1,\dots,h_s)\in (Hilb^{q_1,p_l,\dots,p_1}(X)
\x \Pi_{i=1}^{i=s} 
Hilb^{n_i,p_l,\dots,p_1}(X))$
sont $[h],[h_1],\dots,[h_s]$.
Soit $Z$ le sous-sch\'ema  d\'efini par le produit $I(U_1).\dots. I(U_s)$
des id\'eaux des $U_i$. On pose 
$Inc_{\eta'}=Inc_{ Hilb^{q_1,p_l,\dots,p_1}(X)\x \Pi_{i=1}^{i=s} 
Hilb^{n_i,p_l,\dots,p_1}(X)  } (U,Z)$. 
\\
\textit{Le cas $r>1$.}
Le ferm\'e $Inc^-$ de $Hilb^{q_{r-1},\dots,q_1,p_l,\dots,p_1}(X) \x
\Pi_{i=1}^{ s} Hilb^{n_i,p_l,\dots,p_1}(X)$  d\'efini  
par la  r\'ecurrence se projette sur $\Pi_{i=1}^{ s} Hilb^{n_i,p_l,\dots,p_1}(X)$. On d\'efinit 
$Inc_{\eta'}$ comme le sch\'ema de Hilbert relatif 
$Hilb_{\Pi_{i=1}^{s} Hilb^{n_i,p_l,\dots,p_1}(X)}
^{q_r}(Inc^-)$ de cette projection. 
\\
Si $\eta'\inc \eta$, notons $Inc_{\eta',\eta}$ le ferm\'e de 
$H_{\eta}(X)$ image inverse de $Inc_{\eta'}$ par la projection 
naturelle $H_{\eta}(X)\fd H_{\eta'}(X)$. 
Enfin, posons  $$R_{\eta}(X)=\cap _{\eta'\inc \eta, I(\eta')=1}Inc_{\eta',\eta}.$$

\begin{prop}\label{adherence inclus dans incidence}
Si $X$ est une vari\'et\'e, le plongement 
$\overline{f_{\eta}(F(X,n))} \fd H_{\eta}(X)$ 
se factorise par $R_{\eta}(X)$. 
\end{prop}
\noindent
\textit{D\'emonstration:}
on veut voir que si $\eta'=\{\sigma,\sigma_1,\dots,\sigma_s\}\inc
\eta$ v\'erifie $I(\eta')=1$ et si $x\in F(X,n)$, alors 
$f_{\eta'}(x)\in Inc_{\eta'}$. On proc\`ede 
par r\'ecurrence sur la diff\'erence $d$  entre le niveau de  
$\sigma$ et celui des $\sigma_i$. \\
\textit{Le cas $d=0$}. On a $\sigma=\{\tau_1,\dots, \tau_q\}$,
$\sigma_i=\cup_{k\in I_i} \tau_{k}$. Avec les notations
utilis\'ees dans la d\'efinition de $Inc_{\eta'}$, on a 
$[h]=\cup_{i\leq q}f_{\tau_i}(x)$ et $[h_i]=\cup _{k \in I_i}
f_{\tau_k}(x)$ et on veut 
$I([h])\supset \Pi I([h_i])$.
L'hypoth\`ese $I(\eta')=1$
se traduit par $\cup_i I_i\supset\{1,\dots,q\}$ et implique 
l'inclusion voulue des id\'eaux. 
\\
\textit{Le cas $d>0$}.
Toujours avec les notations
utilis\'ees dans la d\'efinition de $Inc_{\eta'}$,
il nous faut 
v\'erifier que $(f_{\sigma}(x),f_{\sigma_1}(x),\dots,f_{\sigma_s}(x))
\in Hilb^{q_r}(Inc ^-)$. Avec  $\sigma=\{\tau_1,\dots,\tau_q\}$,
cela revient \`a voir pour tout $i$ que $(f_{\tau_i}(x),f_{\sigma_1}(x),
\dots, f_{\sigma_s}(x))\in Inc ^-$. Puisque $I(\eta')=1$, on a pour tout
$i$,
$I(\tau_i,\sigma_1,\dots,\sigma_s)=1$ et il suffit donc d'appliquer 
l'hypoth\`ese de r\'ecurrence. 
\findem

\begin{rem}
Il arrive que l'inclusion $\overline{f_{\eta}(F(X,n))}\inc R_{\eta}(X)$
soit stricte.
\end{rem}
\noindent
\textit{D\'emonstration}: prenons $E=\{1,2,3,4\}$, $\sigma_1=\{1,2\}$,
$\sigma_2=\{1,3\}$, $\sigma_3=\{1,2,3,4\}$ et
$\eta=\{\sigma_1,\sigma_2,
\sigma_3\}$.
Soient $x\in F(X,4)$ et $p=(d,d',q)$ le point de $H_{\eta}(X)=Hilb^2(X)\x
Hilb^2(X)
\x Hilb^4(X)$ d\'efini par $[d]=x_1 \cup x_2,\ [d']=x_1 \cup x_2$ 
et $q=x_1 \cup x_2 \cup x_3 \cup x_4$. Le point $p$ est dans
$R_{\eta}(X)$ mais pas dans $\overline{f_{\eta}(F(X,4))}$ puisque le 
quadruplet est form\'e de quatre points distincts tandis que la r\'eunion 
du support des deux doublets ne contient que deux points. 
\findem
\nl
Dans la suite de ce papier, on s'int\'eressera aux compactifications 
qui sont sympathiques au sens o\`u elles peuvent \^etre d\'efinies 
g\'eom\'etriquement par des conditions d'incidence, ce qui nous
conduit 
\`a la d\'efinition suivante.

\begin{defi}
  On appelle enrichissement admissible un enrichissement 
tel que pour toute vari\'et\'e lisse irr\'eductible,  
$\overline{f_{\eta}(F(X,n))}=  R_{\eta}(X)$.
\end{defi}

\begin{rem}
Nous n'avons en fait pas donn\'e  toutes les relations 
d'incidence possibles. On peut introduire plusieurs fonctions 
d'incidence jouant le r\^ole de la fonction $I$ 
et pour chacune d'elles, construire un lieu d'incidence
correspondant. La d\'efinition g\'en\'erale n\'ecessite un formalisme
lourd. Le point de vue adopt\'e  est d'introduire le formalisme 
minimum qui permette d'atteindre les applications voulues 
(construction d'une classe de compactifications qui contienne les 
vari\'et\'es de collisions et les compactifications classiques). 
\end{rem}

\subsection{D\'efinition par foncteurs repr\'esentables}
\label{subsec:def par foncteurs}
Dans cette section, on explique pourquoi les sch\'emas $R_{\eta}(X)$
peuvent \^etre d\'efinis par un  foncteur repr\'esentable $F_{\eta}$. 
Cette d\'efinition a le double avantage 
d'\^etre extr\^emement maniable et de s'appliquer \`a tous les 
sch\'emas en dehors du cadre des vari\'et\'es r\'eduites. 
\\
Le foncteur $F_{\eta}$ sera d\'efini par  $F_{\eta}(B)=\{$ familles 
$F_1,F_2 \dots$ de sous-sch\'emas 
param\'etr\'ees par $B$ satisfaisant des conditions d'incidence$\}$. 
D\'efinissons maintenant ces conditions d'incidence. 

\begin{prop}\label{def de Inc souligne}
On peut d\'efinir pour tout sch\'ema projectif $P$
et pour tous sous-sch\'emas ferm\'es
$Z\inc B\x Hilb^{q_r,\dots,q_1}(P)$ et $T\inc B\x P$,
tels que $Z$ soit plat sur $B$
un sous-sch\'ema ferm\'e $\underline{Inc}_B(Z,T)$ de $B$
satisfaisant la propri\'et\'e suivante:
un morphisme 
$\phi:S \fd B$ se factorise par  $\underline{Inc}_B(Z,T)$ 
ssi $\underline{Inc}_S(S\x _B Z,S\x _B  T)=S$
\end{prop} 
\noindent
\textit{D\'emonstration:} si $r=0$, $Z$ et $T$ sont des sous-sch\'emas 
de $B\x P$ et il suffit de choisir $\underline{Inc}_B(Z,T)=Inc(Z,T)$ 
en vertu de la proposition \ref{def de Inc}. \\
Si $r\geq 1$, on dispose d'un ferm\'e 
$$\tilde T \inj (B\x Hilb^{q_r,\dots,q_1}(P))\;\x\; P$$
image inverse 
de $T$ par la projection sur $B\x P,$
et d'un ferm\'e  
$\tilde U \inj (B\x Hilb^{q_r,\dots,q_1}(P))\;\x\; Hilb^{q_{r-1},\dots,q_1}(P)$
image inverse de  la famille universelle de 
$Hilb^{q_r,\dots,q_1}(P)\x\; Hilb^{q_{r-1},\dots,q_1}(P) $. 
On peut par r\'ecurrence d\'efinir 
le lieu 
$W=\underline{Inc}_{ B\x Hilb^{q_r,\dots,q_1}(P)}(\tilde U, \tilde T)$.
On d\'efinit alors 
$\underline{Inc}_{B}(Z,T):=
Inc_B (Z,W)$.
\\
Une r\'ecurrence sur $r$ 
montre que si $\phi:S\fd B$ est un changement de base, 
$\underline{Inc}_B(Z,T)\x _B S=\underline{Inc}_S(Z\x _B S, T\x_B S)$.
En particulier, on a $\underline{Inc}(S\x_B Z,S\x_B T)=S 
\Leftrightarrow \underline{Inc}_B(Z,T)\x _B S =S
\Leftrightarrow $ 
$\phi$ se factorise en $S \fd \underline{Inc}_B(Z,T) \fd B$.\findem

\begin{defi}
  On dit que les sous-sch\'emas $Z\inc  B\x Hilb^{q_r,\dots,q_1}(P)$
et  $T_i\inc B\x P$, $i=1,\dots,s$ v\'erifient la relation 
$Z \inc T_1.T_2. \dots .T_s$ si $\underline{Inc}_B(Z,T)=B$, o\`u $T
\inc B\x P$ est d\'efini par l'id\'eal produit $\Pi_{i=1}^{i=s}
I(T_i)$. La relation  $p\inc p_{1}.p_2. \dots . p_{s}$
prend du sens en particulier 
quand $p\in Hilb^{q_r,\dots,q_1,p_l,\dots,p_1}(X)$, $p_i\in 
Hilb^{n_i,p_l,\dots,p_1}(X)$ en prenant $P=Hilb^{p_l,\dots,p_1}(X)$.
\end{defi}

\begin{rem}
\label{rq: Reta est un lieu d'incidence}
Il ressort de cette d\'efinition que $R_{\eta}(X)$ est le sous-sch\'ema 
naturel de $H_{\eta}(X)=H_{\sigma_1}(X)\x \dots \x H_{\sigma_s}(X)$ dont les
points ferm\'es sont les points $(p_1,\dots,p_s)$ 
v\'erifiant $p_{i_0}\inc p_{i_1}. \dots . p_{i_k}$ pour tout 
$\eta'=(\sigma_{i_0},\dots,\sigma_{ i_k})$ tel que $I(\eta')=1$.
\end{rem}

\begin{defi} 
Si $\eta=\{\sigma_1,\dots,\sigma_s\}$ 
est un enrichissement de $E$
avec $\sigma_i\in \Sigma_{n_i,p_{i,s(i)},\dots,p_{i,1}}(E)$,
 on d\'efinit le foncteur  $F_{\eta}$ des sch\'emas sur $k$ dans 
les ensembles par:
  \[F_{\eta}(B)=\{(Z_1,\dots,Z_s) , Z_i \inc B\x H_{\sigma_i}^-(X),
\mbox{ plat et fini de degr\'e $n_i$ sur $B$, } \]
  \[\mbox{ t.q. } I(\sigma_{i_0},\sigma_{i_1},\dots,\sigma_{i_k})=1 \Rightarrow Z_{i_0}
 \inc Z_{i_1}. Z_{i_2}. \dots 
. Z_{i_k}\}\]
\end{defi}

\begin{prop}\label{F_eta representable}
  Le foncteur $F_{\eta}$ est repr\'esentable par le sch\'ema $R_{\eta}(X)$.
\end{prop}
\noindent
\textit{D\'emonstration}: on veut voir qu'il est \'equivalent de se 
donner un \'el\'ement de $F_{\eta}(B)$ ou un morphisme $\phi:B\fd
R_{\eta}(X)$.
Par propri\'et\'e universelle du sch\'ema de
Hilbert,
les ferm\'es $Z_i$ d\'efinissant un \'el\'ement de $F_{\eta}(B)$
correspondent \`a un morphisme $B\fd H_{\eta}(X)$. 
Les relations d'incidence entre les $Z_i$, la propri\'et\'e
universelle 
du sch\'ema $\underline{Inc}$ (prop. \ref{def de Inc souligne}) 
et la remarque 
\ref{rq: Reta est un lieu d'incidence}  montrent que ce morphisme 
se factorise en un morphisme $\phi:B\fd R_{\eta}(X)$. R\'eciproquement, 
il est clair par les m\^emes arguments qu'un morphisme $\phi:B\fd R_{\eta}(X)$
d\'efinit des $Z_i$ satisfaisant les relations d'incidence.
\findem

\section{Propri\'et\'es des compactifications r\'esultant des
  d\'efinitions}
\label{sec:prop resultant des def}
Dans cette section, on exploite la d\'efinition fonctorielle des 
sch\'emas $R_{\eta}(X)$ pour en donner les premi\`eres propri\'et\'es. 
Pour chacune des propri\'et\'es, on explique comment elle sera 
utilis\'ee lors du  th\'eor\`eme de classification.

\subsection{Changements de bases }
Si $\eta$ et $\eta'$ sont deux enrichissements,
pour montrer que pour tout $X$, $R_{\eta}(X)=R_{\eta'}(X)$,
il suffit de le montrer pour $X$ lisse irr\'eductible: 

\bprop
\label{changement de base}
Soit $Y$ un sous-sch\'ema localement ferm\'e d'un sch\'ema $X$. Les 
sch\'emas $R_{\eta}(Y)$ et 
$R_{\eta}(X) \x _{Hilb^n(X)} Hilb^n(Y)$ sont canoniquement 
isomorphes.
\eprop
v\'erifions que ces deux sch\'emas repr\'esentent le m\^eme foncteur. 
Notons $F$ et $G$ les foncteurs associ\'es \`a 
$$R_{\eta}(Y) \mbox{ et } R_{\eta}(X) \x _{Hilb^n(X)} Hilb^n(Y)$$
Soit $B$ un $k$-sch\'ema. On veut \'etablir une bijection
 canonique entre $F(B)$ et $G(B)$. Posons
 $\eta=\{\sigma_1,\dots,\sigma_s\}$, o\`u $\sigma_s=\{1,\dots,n\}$.
Soit $(Z_1,\dots,Z_s)\in F(B)$. Alors  $(Z_1,\dots,Z_s,Z_s)\in G(B)$.
R\'eciproquement si $(Z_1,\dots,Z_s,Z)\in G(B)$, alors $Z_s=Z\subset 
B\x Hilb^n(Y)$ puisque $Z$ et $Z_s$ d\'efinissent le m\^eme morphisme
$B \fd H_{\eta}(X)$. Et puisque $Z_i \inc Z_s$, $Z_i \inc B\x
H_{\sigma_i}^-(X)$ est en fait dans $B\x H_{\sigma_i}^-(Y)$ 
et $(Z_1,\dots,Z_s,Z_s)\in F(B)$.
L'identification de $F(B)$ et $G(B)$ en d\'ecoule.
\findem

\subsection{Morphismes d'oubli}

\bprop \label{existence du morphisme oubli}
Soient $\eta$ et $\eta'$ deux enrichissements avec $\eta \inc \eta'$. 
Il existe un morphisme d'oubli $p_{\eta' ,\eta}:R_{\eta'}(X)
\fd R_{\eta}(X)$. De plus, $p_{\eta',\eta}\circ f_{\eta'} =f_{\eta}$.
\eprop
la projection naturelle  $p_{\eta',\eta}$ 
de $H_{\eta'}(X)$ sur $H_{\eta}(X)$
envoie $R_{\eta'}(X)$ dans $R_{\eta}(X)$ par d\'efinition
de ces lieux d'incidence.  L'egalit\'e  
$p_{\eta',\eta}\circ f_{\eta'} =f_{\eta}$ est  \'evidente.
\findem \nl
La proposition suivante nous dit que, sous certaines conditions,
le morphisme d'oubli est un isomorphisme. Elle r\'eduira drastiquement 
le nombre de $R_{\eta}(X)$ \`a \'etudier lors de la classification. 

\bprop
\label{oubli=isomorphisme residuel}
Soit $\eta$ un enrichissement admissible et $\eta' =\eta \cup {\sigma}$. 
Supposons que $\eta$ contienne deux enrichissements
$\sigma_1$ et $\sigma_2$ de m\^eme niveau que $\sigma$ v\'erifiant
\beq
\sigma &\in& \Sigma_{1,p_{l-1},p_{l-2}, \dots,p_1}(E)
=\Sigma_{p_{l-1},p_{l-2}, \dots,p_1}(E)\\
\sigma_1 &\in& \Sigma_{p_l,p_{l-1},,p_{l-2}, \dots,p_1}(E)\\
\sigma_2 &\in& \Sigma_{p_l-1,p_{l-1},,p_{l-2}, \dots,p_1}(E)\\
\sigma _1 &=&\sigma_2  \cup \{\sigma\}
\eeq
alors pour tout $X$, le morphisme d'oubli 
$$p_{\eta' ,\eta}:R_{\eta'}(X) 
\fd R_{\eta}(X)$$
est un isomorphisme.
\end{prop}
\blm [existence d'une application r\'esiduelle]\ \\
\label{residuelenfamille}
Soient $F_1$ et $F_2$ deux familles  de $B\x X$ plates et finies 
sur une base $B$ v\'erifiant $F_1 \inc F_2$. Supposons
que les polyn\^omes de Hilbert des fibres de  $F_1$ et $F_2$ soient 
respectivement $p$ et $p+1$. Il existe un unique morphisme 
$res(F_1,F_2):B \fd X$ dont le graphe $Res(F_1,F_2)$ v\'erifie 
$$ I(F_1).I(Res(F_1,F_2))\inc I(F_2)$$
On dira que $res(F_1,F_2)$ est l'application r\'esiduelle d\'efinie par l'inclusion
$F_1 \inc F_2$. En particulier, pour deux sous-sch\'emas de $X$ 
de dimension 
z\'ero $Y$ et $Y'$ v\'erifiant $Y\inc Y'$ et 
$col(Y')=col(Y)+1$, il existe un sch\'ema 
$Res(Y,Y')$.
\end{lm} \noindent
\textit{D\'emonstration du lemme}:
montrons que le faisceau d'id\'eaux $(I(F_2):I(F_1))$ d\'efinit  un
ferm\'e $Z$ de $B\x X$, plat 
sur $B$, dont la 
fibre au dessus de tout point $b$ de $ B$ est 
de degr\'e un. Le probl\`eme est local sur $B$ qu'on peut 
supposer affine. Les sch\'emas $F_1$ et $F_2$ sont alors
\'egalement affines. Puisqu'on a l'inclusion 
$$(I(F_2):I(F_1))
\supset I(F_2)$$
le sous-sch\'ema $Z$ d\'efini par 
$(I(F_2):I(F_1))$ est un sous-sch\'ema ferm\'e 
de $F_2$. L'id\'eal de ${\cal O}_{B\x X}/I(F_2)$ d\'efinissant 
$Z$ comme sous-sch\'ema de $F_2$ est $ (I(F_2):I(F_1))/I(F_2)
=(0:(I(F_1)/I(F_2)))$. De la suite exacte 
$$0\fd I(F_1)/I(F_2) 
{\fd} \Gamma(F_2)
{\fd} \Gamma(F_1)\fd 0$$
on d\'eduit que $I(F_1)/I(F_2)$ est un $\Gamma(B)$-module
localement libre de rang un et, quitte \`a restreindre $B$, on 
peut supposer qu'il est  libre engendr\'e par un 
\'el\'ement $f$ de $\Gamma(F_2)$. Le diagramme commutatif 
\`a ligne exacte 
\beq
\barr{cccccc}
&&\Gamma(Z)&&\\
& \fhd&&\fbd&\\
0\fd \;\;0:(I(F_1)/I(F_2))\;\;
\fd \Gamma(F_2)&& \stackrel{\times f}
{\fd}&&\Gamma(F_2)
\earr
\eeq
montre que  $\Gamma(Z)=\Gamma(F_2)/(0:I(F_1)/I(F_2))$ s'identifie au  sous $\Gamma(B)$-module de 
$\Gamma(F_2)$ form\'e par les multiples de $f$.
Finalement la suite exacte
$$0\fd \Gamma(Z) \fd \Gamma(F_2)\fd \Gamma(F_1)\fd 0$$
montre que $\Gamma(Z)$ est un $\Gamma(B)$-module plat, donc localement libre,
de rang un.
L'inclusion de $Z$ dans $B\x X$ induit  un 
morphisme de $B$ dans $X$ ayant la 
propri\'et\'e voulue.
\nl 
Soit $res(F_1,F_2)$ un morphisme ayant la propri\'et\'e 
voulue. On a l'inclusion  $I(Res(F_1,F_2)) \inc (I(F_2):I(F_1))$.
Deux familles, plates et finies de degr\'e un sur la base et incluses l'une dans l'autre sont n\'ecessairement
\'egales donc $I(Res(F_1,F_2))= (I(F_2):I(F_1))$.
 \findem
\brq
Dans [LB], Le Barz a montr\'e l'existence d'un r\'esiduel
pour deux sous-sch\'emas $Y$ et $Y'$ de dimension z\'ero de $X$. Le 
lemme est une version relative du r\'esultat 
de Le Barz
\erq
\textit{D\'emonstration de la proposition 
\ref{oubli=isomorphisme residuel}}: quitte \`a r\'ealiser $X$ comme sous-sch\'ema
d'une vari\'et\'e lisse irr\'eductible et \`a utiliser la proposition 
\ref{changement de base}, on peut supposer $X$ lisse irr\'eductible. Posons $\eta=\{\sigma_1,\sigma_2,\dots,
\sigma_s\}$.
Le sch\'ema $R_{\eta}(X)$ est muni de ferm\'es universels 
$U_1,U_2,\dots,U_{s}$. On d\'efinit le ferm\'e $U_{\sigma}$ de 
$R_{\eta}(X)\x H_{\sigma}^-(X)$
par l'id\'eal $I(U_{\sigma}):=(I(U_1):I(U_2))$.
Les ferm\'es 
universels 
$U_1, U_2, \dots, U_s, \ U_{\sigma}$ d\'efinissent 
par  propri\'et\'e universelle  
 un morphisme 
$p_{\eta,\eta'}:R_{\eta}(X)\fd R_{\eta'}(X)$. 
Les morphismes 
$p_{\eta,\eta'}$ et $p_{\eta',\eta}$ sont par construction inverses 
l'un de l'autre.\findem

\subsection{Actions de groupes}

On dispose d'une action naturelle du groupe sym\'etrique $S_n$ sur l'ensemble
$\{1,\dots,n\}$. Cette action induit une action de $S_n$ sur l'ensemble 
des structures de $\{1,\dots,n\}$ par permutation des indices, puis une 
action sur l'ensemble des enrichissements par permutation des structures.
On a la proposition \'evidente suivante.
\bprop
\label{isomorphisme par permutation}
Soient $\eta$ et $\eta'$ deux enrichissements de $\{1,\dots,n\}$ qui sont 
dans la m\^eme orbite pour l'action de $S_n$. Alors les sch\'emas 
$R_{\eta}(X)$ et $R_{\eta'}(X)$ sont isomorphes.
\end{prop}
\noindent
Si $\eta=\{\sigma_1,\dots,\sigma_r\}$, on note $G_{\eta}:=\{g \in S_n 
\mbox{ t.q. }\exists p \in S_r,\ g.\sigma_i=\sigma_{p(i)}\}$ et 
$H_{\eta}:=\{g \in S_n \mbox{ t.q. }\forall i, \ g.\sigma_i \linebreak[0]
= \sigma_i\}$.

\bprop \label{description du groupe agissant}
Pour tout 
$X$, il existe une action naturelle du quotient  
$G_{\eta}/H_{\eta}$ sur le sch\'ema $R_{\eta}(X)$
\eprop
si $\eta=\{\sigma_1,\dots,\sigma_r\}$, un morphisme $Z \fd R_{\eta}(X)$ 
est d\'efini par la donn\'ee de ferm\'es $U_i \inc 
Z \x H_{\sigma_i}^{-}(X)$ ($1 \leq i \leq r$)  
satisfaisant \`a des relations d'incidence. En particulier le
morphisme identit\'e de $R_{\eta}(X)$ d\'efinit
un ensemble de ferm\'es
$U_1,\dots,U_r$. Toute permutation  $g$
de $G_{\eta}$ est associ\'e \`a une permutation $p$.
L'ensemble  $U_{p(1)},\dots, 
U_{p(r)} $  d\'efinit par propri\'et\'e universelle 
un endomorphisme $\phi_g$ de $R_{\eta}(X)$, qui est 
un automorphisme 
d'inverse $\phi_{g^{-1}}$. L'application 
$g \fd \phi_g$ d\'efinit une action de groupe de $G_{\eta}$ sur 
$R_{\eta}(X)$. Puisque $H_{\eta}\inc G_{\eta}$ n'agit pas,
l'action se factorise en une action
de $G_{\eta}/H_{\eta}$.
\findem

\section{Etude et classification des compactifications pour $n \leq 3$.}
\label{sec:classification}
On veut classer les compactifications \`a isomorphisme pr\`es 
(au sens o\`u 
deux compactifications  $C$ et $C'$ d'un m\^eme 
quotient $F(X,n)/Q$ sont isomorphes s'il existe 
un isomorphisme entre $C$ et $C'$ valant l'identit\'e sur 
$F(X,n)/Q$). En d'autres termes, on cherche une liste d'enrichissements
admissibles tels que pour tout $\eta$ admissible, $R_{\eta}(X)$ est
isomorphe 
\`a $R_{\beta}(X)$ pour un unique $\beta$ dans la liste des admissibles. 
Cette section donne la classification pour $n=2$ 
et $n=3$.

\subsection{Le cas $n=2$.}
Le cas $n=2$ \'etant facile, nous nous contentons de 
donner le r\'esultat. Les d\'emonstrations sont des 
cas particuliers extr\^emement simples du cas $n=3$ 
et nous les omettons par souci de concision. 
\nl
Notons $\eta_0$ l'enrichissement de $\{1,2\}$ 
contenant l'unique structure $\{1,2\}$,
et $\eta_1$   l'enrichissement contenant les deux structures $\{1,2\}$
et $\{1\}$.

\begin{thm}
 Tous les enrichissements de $\{1,2\}$ sont admissibles. 
Quel que soit $X$ et quel soit l'enrichissement $\eta$, $R_{\eta}(X)$
est isomorphe \`a $R_{\eta_0}(X)$ ou \`a $R_{\eta_1}(X)$.
De plus, $R_{\eta_0}(X)=Hilb^2(X)$ et, pour $X$ lisse,
$R_{\eta_1}(X)$ est l'\'eclatement 
de $X\x X$ le long de la diagonale. 
\end{thm}

\subsection{Le cas $n=3$}
La classification s'effectue en quatre \'etapes
\ben
\item exhiber ``\`a la main '' quelques enrichissements admissibles
\item  exhiber quelques enrichissements non admissibles
\item \'etablir des lemmes de contamination qui permettent, \'etant 
donn\'e un enrichissement admissible $\eta$, de montrer que 
d'autres enrichissements $\eta'$ sont admissibles.
\item traiter tous les enrichissements \`a partir des cas particuliers 
et des lemmes de contamination.
\een

\subsubsection{Notations}
Dans cette section, on met en place les notations qui nous facilitent 
la manipulation des enrichissements. 
Notons $E=\{1,2,3\}$. On note les structures de niveau un de $E$ \`a l'aide
d'indices ($\{i,j,k\}=\{1,2,3\}$):
$$\sigma_i=\{i\} \in \Sigma_1(E)$$
$$\sigma_{ij}=\{i,j\} \in \Sigma_2(E)$$
$$\sigma_{123}=\{1,2,3\} \in \Sigma_3(E)$$
On note les structures de niveau deux \`a l'aide d'exposants:
$$\sigma^k=\{\sigma_{ik},\sigma_{jk}\} \in \Sigma_{2,2}(E)$$
$$\sigma^{123}=\{\sigma_{12,},\sigma_{13},\sigma_{23}\} \in \Sigma_{3,2}(E)$$
On repr\'esente les enrichissements avec des indices et des exposants, 
suivant les structures contenues dans l'enrichissement
$$\eta^{a_1,a_2,\dots,a_r}_{b_1,b_2,\dots,b_s}:=\{\sigma^{a_1},\sigma^{a_2},
\dots,\sigma^{a_r},\sigma_{b_1},\sigma_{b_2}, \dots,\sigma_{b_s}\}$$
On utilisera \'egalement cette notation avec indices et exposants 
pour les sch\'emas $\re$ et $H_{\eta}(X)$:
$$R^{a_1,a_2,\dots,a_r}_{b_1,b_2,\dots,b_s}(X):=R_{\eta^{a_1,a_2,\dots,a_r}_{b_1,b_2,\dots,b_s}}(X)$$ 
$$H^{a_1,a_2,\dots,a_r}_{b_1,b_2,\dots,b_s}(X):=H_{\eta^{a_1,a_2,\dots,a_r}_{b_1,b_2,\dots,b_s}}(X)$$ 
Parmi les enrichissements de niveau inf\'erieur ou \'egal \`a deux, il y 
en a un maximum, celui qui contient toutes les structures de niveau 
inf\'erieur ou \'egal \`a  deux. On le notera $\eta_{max}$ et $R_{max}(X)$ le 
sch\'ema de triplets de $X$ correspondant. Enfin, si $p$ est un point
de $X$, on notera $[[p]]$ le voisinage formel de $p$ dans $X$, 
c'est \`a dire le morphisme naturel $\s \hat A \fd X$, o\`u $\hat A$ est la
limite projective des $A/{\goth m}^n$, $(A,{\goth m})$ 
\'etant l'anneau local de
$X$ en $p$.

\subsubsection{Quelques enrichissements admissibles}
Dans cette section, on montre que les enrichissements $\eta_{123},
\eta_{1,2,3,12,13,23,123},$ 
$\eta_{12,123}, \eta ^1_{123}, \eta_{123}^{123}$ sont
admissibles.\\
On sait que pour $X$ lisse irr\'eductible, $Hilb^3(X)$ est irr\'eductible et 
est l'adh\'erence de $F(X,3)/S_3$ donc  $\eta_{123}$
est admissible.\\
Le cas $\eta_{1,2,3,12,13,23,123}$ a \'et\'e trait\'e dans [LB],  
o\`u Le Barz a montr\'e 
que cet enrichissement \'etait admissible.  
\\
Pour les trois enrichissements $\eta$ restant, 
on montre que si $X$ est lisse
irr\'eductible,
les $R_{\eta}(X)$ sont irr\'eductibles de dimension $3.dim(X)$.
Cela suffit \`a montrer que les enrichissements sont admissibles
car l'unique composante irr\'eductible de $R_{\eta}(X)$ est 
n\'ecessairement la 
compactification de $f_{\eta}(F(X,3))$. Pour 
obtenir l'irr\'eductibilit\'e,
on montre que les  $R_{\eta}(X)$  sont lisses connexes. Ils admettent 
des points sp\'eciaux et la lissit\'e se montre par un calcul 
en coordonn\'ees locales au voisinage de ces points sp\'eciaux.
Le dimension de la composante irr\'eductible est une cons\'equence de 
ce calcul local. 
\bprop
\label{point special}
Soit $\eta$ un enrichissement de $\{1,2,3\}$ de niveau au plus deux. 
Soit $X$ une vari\'et\'e lisse de dimension au moins deux.
Soient $x$ un point de $X$, $[d]$ un sous-sch\'ema ponctuel de $X$ de 
colongueur
deux contenant $x$. Soit $[t]$ un sous-sch\'ema ponctuel de 
colongueur trois 
de $X$ contenant $[d]$ isomorphe \`a $\s k[x,y]/(x^2,xy, y^2)$ en tant que 
sch\'ema abstrait. Il existe un point $q(x,d,t)$ de $R_{\eta}(X)$ 
enti\`erement d\'etermin\'e par $x,d$ et $t$. 
\eprop
nous allons exhiber un 
point de $R_{\eta_{max}}(X)$ d\'etermin\'e par $x,d$ et $t$. 
L'image de ce point par le morphisme d'oubli $p_{\etamax,\eta}$ 
sera le point 
$q(x,d,t)$. 
Les doublets inclus dans le  triplet $[t]$ forment un sous-sch\'ema 
de $\hdx$ isomorphe \`a  $\PP^1$, et 
donc un doublet de doublets 
inclus dans le triplet correspond \`a un sous-sch\'ema de degr\'e deux du $\PP^1$. Le doublet $d$ inclus dans $[t]$ est 
 un point du  $\PP^1$. Il existe un unique sous-sch\'ema de degr\'e 
deux de $\PP^1$ support\'e par $d$. On note $d^2$ le point de $Hilb^{2,2}(X)$
associ\'e. On a par construction $[d^2] \inc [t]$. De m\^eme, on peut d\'efinir un point $d^3$ de $Hilb^{3,2}(X)$
pour lequel $[d^3] \inc [t]$. Consid\'erons alors le point 
$$p=(x,x,x,d,d,d,d^2,d^2,d^2,d^3,t)$$
de 
\beq
  \hachun \x \hachdeux &\x& \hachtrois \x \hachundeux
\x \hachuntrois \x \hachdeuxtrois\\
& \x& H^1(X) \x H^2(X) \x H^3(X) \x H^{123}(X) \x H_{123}(X) 
\eeq
On v\'erifie que toutes 
les relations d'incidence sont satisfaites et donc $p \in R_{\eta_{max}}(X)$.
\findem

\bprop
\label{connexite}
Soient $X$ une vari\'et\'e lisse connexe et $\eta$ un enrichissement de 
$\{1,2,3\}$ de niveau au plus deux. Le sch\'ema $\re$ est connexe. 
\eprop
soient $q$ un point de $\re$ et $q(x,d,t)$ un point sp\'ecial de $R_{\eta}(X)$ au sens de la proposition
\ref{point special}. Nous allons d\'eterminer une suite de points 
$q_0=q,q_1,\dots,q_n=q(x,d,t)$ telle que, pour tout $i$, il existe une 
courbe $C_i$ dans $\re$ contenant $q_i$ et $q_{i+1}$.\\
\newcommand{\emx}{\eta_{\max}}
Si $\eta=\emx$, le point $q$   est  d\'efini 
par ses coordonn\'ees 
$$p_1(q),p_2(q),p_3(q),p_{12}(q),p_{13}(q),
p_{23}(q),p^1(q),p^2(q),p^3(q),p^{123}(q),p_{123}(q)$$
dans
$$X^3 \x \hdx^3 \x \hdhdx^3 \x \htrhdx \x \htrx$$ 
Pour $\eta$ quelconque, $q$ est d\'etermin\'e par un certain nombre 
de coordonn\'ees  parmi $p_1(q),\dots,p_{123}(q)$.
\\
Par d\'efinition d'un enrichissement, la coordonn\'ee $p_{123}$ fait toujours
partie des coordonn\'ees d\'efinissant $q$.
Le support de $p_{123}(q)$ est constitu\'e d'au plus trois points. En bougeant 
$p_{123}(q)$ le long d'une famille \`a un param\`etre, on peut construire une 
courbe $C_0$ dans $R_{\eta}(X)$ telle que:
$C_0$ contienne $q_0$ et un point $q_1$ dont le triplet $[p_{123}(q_1)]$ 
associ\'e est ponctuel de support $x$. 
\\
Soit $[[x]]$ le voisinage formel de $x$ dans $X$. Le triplet
$[p_{123}(q_1)]$ est maintenant un sous-sch\'ema de $[[x]]$. Puisqu'on a les 
relations d'incidence:
\beq
p_i(q_1) \inc [p_{123}(q_1)]  \inc  [[x]]  \\
\ [ p_{ij}(q_1)] \inc [p_{123}(q_1)] \inc [[x]] \\
\ [p^i(q_1)] \inc [p_{123}(q_1)] \inc [[x]] \\
\ [p^{123}(q_1)] \inc [p_{123}(q_1)] \inc [[x]]
\eeq 
``tout se passe'' dans le 
voisinage formel $[[x]]$. En particulier tout automorphisme 
$\phi_t$ de $[[x]]$ d\'efinit un point $\psi_t(q_1)$ de $\re$.
En tant que sch\'ema abstrait, $[p_{123}(q_1)]$ est soit isomorphe \`a 
$\s k[y,z]/(y^2,yz,z^2)$, 
soit isomorphe \`a $\s k[z]/z^3$. Dans le premier cas, 
on d\'efinit $q_2=q_1$. Dans le deuxi\`eme cas, et si on suppose 
que $X$ est de dimension sup\'erieure ou \'egale \`a trois, 
on peut choisir un 
isomorphisme 
$$[[x]] \simeq \s k[[y,z,t_1,t_2,\dots,t_n]]$$ 
dans lequel $[p_{123}(q_1)]$ ait pour \'equations $(y,z^3,t_1,\dots,t_n)$.
Pour tout $t$ dans $\A^1-0$, on d\'efinit l'automorphisme $\phi_t$ de $[[x]]$:
\beq
y&\mapsto &ty+z^2\\
z &\mapsto & z \\
t_i &\mapsto & t_i
\eeq
On en d\'eduit un morphisme 
\beq
\psi:\AA^1-0 & \fd & \re\\
t & \mapsto & \phi_t(q_1)
\eeq
Il existe un morphisme 
$\tilde{\psi}$ prolongeant $\psi$ en $0$ et $\infty$ et $\tilde{\psi}(\infty)
=q_1$. Notons $q_2=\tilde{\psi}(0)$. On 
v\'erifie par un petit calcul que $[p_{123}(q_2)]$ est d\'efini par 
l'id\'eal $(y^2,yz,z^2,t_1,\dots,t_n)$.
\\
Les doublets inclus dans $[p_{123}(q_2)]$ forment un $\PP^1$ et 
le morphisme de $\PP^1$ dans $\hdx$ correspondant peut \^etre d\'ecrit 
de la fa\c con suivante:
\beq
 \PP^1 &\fd& \hdx\\
(h_0:h_1) &\mapsto& x \mbox{ avec } I([x])=(h_0y-h_1z,y^2,yz,z^2,t_1,t_2,\dots,t_n)
\eeq
Choisissons la famille d'automorphismes suivantes de $[[x]]$. 
\beq 
\phi_t:&y & \mapsto ty\\
&z &\mapsto z\\
&t_i &\mapsto t_i
\eeq
Comme pr\'ec\'edemment, $\phi_t$ induit un morphisme $\tilde{\psi}:\AA^1
\fd \re$. On a $q_2=\tilde{\psi}(1)$ et on cherche \`a comprendre $q_3:=\tilde{\psi}(0)$. Commen\c cons 
par d\'eterminer $p^1(\tilde {\psi}(0))$.
Le point $p^1(\tilde {\psi}(1))$
correspond \`a un sous-sch\'ema de degr\'e deux $D$  du $\PP^1$ form\'e par 
les doublets inclus dans le triplet. 
Le morphisme 
$\phi_t$ induit  un automorphisme $\rho_t$ de $\PP^1$: $$\rho_t(h_0:h_1)=(th_0:h_1)$$
et 
$p^1(\tilde {\psi}(t))$ correspond au sous-sch\'ema $\rho_t(D)$. 
Si on a choisi un bon isomorphisme 
$$[[x]] \simeq \s k[[y,z,t_1,t_2,\dots,t_n]]$$ 
le support de $D$ ne rencontre pas le point \`a l'infini $(1:0)$, 
le doublet $d$ correspond au point $(0:1)$ de $\PP^1$,
et alors $[p^1(\tilde {\psi}(0))]$ est le sous-sch\'ema de degr\'e deux 
de $\PP^1$ support\'e par  $d=(0:1)$.
\\
La m\^eme d\'emonstration montre que tous les points 
 $$p^i(\tilde {\psi}(0)),\
 p_{ij}(\tilde {\psi}(0)),\ p^{123}(\tilde {\psi}(0))$$ correspondent \`a 
des sous-sch\'emas de $\PP^1$ support\'es par $d=(0:1)$, ce qui signifie 
$\tilde {\psi}(0)=q(x,d,t)$.
\findem
\brq
on a fait la d\'emonstration dans le cas o\`u la dimension de $X$ est 
sup\'erieure ou \'egale \`a trois. Dans le cas o\`u elle vaut deux, il suffit
de supprimer les lignes contenant des ``$t_i$'' dans la d\'emonstration. 
Le cas de la dimension  un est facile et laiss\'e au lecteur.
\erq

\bcor \label{critere de lissite}
Si $\re$ admet un point singulier, alors les points $q(x,d,t)$ sont 
singuliers. 
\ecor
\demo: supposons $R_{\eta}(X)$ singulier. La construction pr\'ec\'edente
d\'efinit une famille de courbes $C_i$ joignant $q$ \`a $q(x,d,t)$. 
Pour un point $p$ de $C_0$ diff\'erent de $q_1$, les deux voisinages 
formels $[[p]]$ et $[[q]]$ sont isomorphes. Donc si $q_0:=q$ est
singulier, tous les points de $C_0$ diff\'erents de $q_1$ sont 
singuliers.  Par suite $q_1$ est singulier. Le m\^eme raisonnement 
permet d'obtenir de proche en proche
$q_2,\dots, q_n=q(x,d,t)$ singuliers.
\findem

\begin{coro} \label{critere d'isomorphisme d'oubli}
  Soient  $\eta$ est un enrichissement admissible et $\eta_{res}\inc \eta$
 l'ensemble des structures de niveau au plus deux de $\eta$. Si 
 la fibre d'un morphisme d'oubli $p_{\eta',\eta}:R_{\eta'}(X)\fd R_{\eta}(X)$ est 
sch\'ematiquement r\'eduite \`a un point au dessus des points 
$q$ tels que $p_{\eta,\eta_{res}}(q)$ est de la forme $q(x,d,t)$,
alors $p_{\eta',\eta}$ est un isomorphisme. 
\end{coro}
\demo: la fibre au dessus d'un point g\'en\'eral 
$p=f_{\eta}(x),\ x \in F(X,3)$ est non vide car elle contient
$f_{\eta'}(x)$. Pour montrer l'isomorphisme, il suffit donc 
par semi-continuit\'e de voir 
que toutes les fibres sont incluses dans un point. Toujours par
semi-continuit\'e et en raisonnant comme dans la proposition
pr\'ec\'edente, il suffit de le v\'erifier aux points sp\'eciaux $q$
tels que $p_{\eta,\eta_{res}}(q)$ est de la forme 
$q(x,d,t)$, ce qui est vrai par hypoth\`ese. 
\findem

\bprop \label{R_{12,123}admissible}
Soit $X$ une vari\'et\'e lisse irr\'eductible. Le sch\'ema $R_{12,123}(X)$ 
est une vari\'et\'e lisse irr\'eductible de dimension $3.dim(X)$.
En particulier $\eta_{12,123}$ est admissible.
\eprop
puisque $R_{12,123}(X)$ est connexe (proposition \ref{connexite}), il suffit 
de v\'erifier la lissit\'e,  en les points sp\'eciaux  $q(p,d,t)$ par 
\ref{critere de lissite}.
Supposons la vari\'et\'e $X$ de dimension au moins trois. Pour un tel 
point $q$, on peut choisir un voisinage formel $\s k[[x,y,z_1,z_2,\dots,z_n]]$ 
de $p$ 
dans lequel les \'equations de $[d]$ et $[t]$ sont respectivement:
$$I([d])=(x^2,y,z_1,\dots,z_n)$$
$$I([t])=(x^2,xy,y^2,z_1,z_2,\dots,z_n)$$
D'apr\`es [LB], le voisinage formel de $d$ dans $\hdx$ 
est isomorphe \`a $$\s [[a,b,c,d,e_1,e_2,\dots,e_n,f_1,f_2,\dots,f_n]]$$
et le voisinage formel de $t$ dans $\htrx$ est isomorphe 
\`a $$\s k[[u,u',u'',v,v',v'',\rho_1,\rho_2,\dots,\rho_n,\sigma_1,\sigma_2,
\dots,\sigma_n,\theta_1,\theta_2,\dots,\theta_n]]$$
Les id\'eaux universels de $\hdx$ et $\htrx$ au dessus de ces voisinages 
formels sont:
$$I_{12}=(x^2+ax+b,y-cx-d,z_i-e_ix-f_i)$$
$$I_{123}=(x^2+ux+vy+w,xy+u'x+v'y+w',y^2+u''x+v''y+w'',z_i+\rho_i x +\sigma_i y+\theta_i)$$
o\`u $w,w',w''$ sont des fonctions alg\'ebriques de $u,u',u'',v,v',v''$.
Le voisinage formel de $(d,t)$ dans $R_{12,123}(X)$ est le lieu  
 $$ Z \inj \s k[[a,b,c,d,e_i,f_i,u,u',u'',v,v',v'',\rho_i,\sigma_i,\theta_i]]$$
au dessus duquel $I_{123} \inc I_{12}$.\\
Le calcul de $Z$ a d\'ej\`a \'et\'e effectu\'e dans $[LB]$ (Le Barz y 
parlait du lieu ensembliste mais il a effectu\'e le calcul au moyen 
de divisions et les techniques d\'evelopp\'ees dans [Ev] 
montrent qu'il a en fait calcul\'e les lieux sch\'ematiques 
d'incidence). L'id\'eal de $Z$ est le suivant:
\beq
I(Z)=(u-a+cv,b-dv-w,u'-ac+cv'+d,2cd+cv''+u''-ac^2,\\
e_i+\sigma_ic+\rho_i,f_i+\theta_i+\sigma_id)
\eeq
et
\beq
Z&\simeq&  \s [[a,b,c,d,e_i,f_i,u,u',u'',v,v',v'',\rho_i,\sigma_i,\theta_i]]/I(Z)\\
&\simeq &\s k[[a,c,d,v,v',v'',\rho_i,\sigma_i,\theta_i]]
\eeq
Donc $Z$ est lisse au point sp\'ecial et a la dimension attendue.
\nl
On a fait la d\'emonstration dans le cas o\`u $X$ est de dimension au moins 
trois. Dans le cas o\`u $X$ est de dimension deux, il suffit d'\^oter de 
la d\'emonstration les lignes contenant des termes $z_i$. Le cas de la 
dimension un est facile et laiss\'e au lecteur.
\findem

\bprop
\label{R_{123}^1 admissible}
Si $X$ est une vari\'et\'e lisse irr\'eductible, le sch\'ema $R_{123}^1(X)$ 
est une vari\'et\'e lisse irr\'eductible de dimension $3.dim(X)$. En
particulier $\eta_{123}^1$ est admissible.
\eprop
pour les m\^emes raisons que 
pr\'ec\'edemment, nous allons montrer la lissit\'e au point sp\'ecial
dans le cas o\`u $X$ est de dimension deux.
\nl
Soit $(p^1,p_{123})$ un point sp\'ecial de $R_{123}^1(X)$.
Les objets en jeu sont un sch\'ema ponctuel $[p_{123}]$ de $X$ de colongueur 
trois, un sch\'ema $[p^1]$ de $\hdx$ ponctuel de colongueur deux 
inclus dans $[p_{123}]$, un sch\'ema 
$[p_{12}]$ ponctuel de $X$ o\`u $p_{12}$ est le support de $[p^1]$,
un point $p$ support \`a la fois de $[p_{12}]$ et de  $[p_{123}]$.
\nl
On peut choisir un voisinage formel $\s k[[x,y]]$ de $p$ tel que:
\bit
\item
$I([p_{12}])=(x^2,y)$, 
le voisinage formel de $p_{12}$ dans $Hilb^2(X)$ est isomorphe \`a 
$$
\s k[[a,b,c,d]]
$$
et 
l'id\'eal universel de $\hdx\x X$ au dessus de cette carte est 
$$(x^2+ax+b,y-cx-d)$$
\item
$I([p_{123}])=(x^2,xy,y^2)$, 
le voisinage formel de $p_{123}$ est isomorphe \`a 
$$\s k[[u,u',u'',v,v',v'']]$$
et l'id\'eal universel de $\htrx \x X$ au dessus de cette carte est
$$I_{123}=(x^2+ux+vy+w,xy+u'x+v'y+w',y^2+u''x+v''y+w'')$$ o\`u $w,w',w''$ sont des fonctions alg\'ebriques de $u,u',u'',v,v',v''$.
\item
le voisinage formel de $p^1$ dans $Hilb^{2,2}(X)$ est isomorphe \`a 
$$
\s k[[e,f,g,h,i,j,l,m]]
$$ 
et l'id\'eal universel de $Hilb^{2,2}(X) \x \hdx $ au dessus
de cette carte est:
$$I^1=(a-ec-f,b-gc-h,c^2-ic-j,d-lc-m)$$
car $[p^1]$ a pour \'equations $a=b=c^2=d=0$ dans $Hilb^2(X)$
\eit
D'apr\`es la d\'emonstration pr\'ec\'edente, le lieu $Z$ de $\hdx \x \htrx$ 
au dessus duquel $I_{123} \inc I_{12}$  est donn\'e par l'id\'eal 
$$J=(u-a+cv,b-dv-w,u'-ac+cv'+d,2cd+cv''+u''-ac^2)$$
En prenant l'image inverse par les projections \'evidentes, $I^1$ et 
$J$ peuvent \^etre vus comme des id\'eaux de $(Hilb^{2,2}(X) \x \htrx)\ \x \ \hdx$.
Par d\'efinition du lieu d'incidence, le lieu $W$ de $Hilb^{2,2}(X) \x \htrx$ au dessus
duquel $[p^1] \inc [p_{123}]$ est le lieu au dessus duquel $J \inc I^1$.
\\
Gr\^ace au morphisme canonique $\s k[[a,b,c,d]] \fd \s k[a,b,c,d]$, on peut 
supposer que $I^1$ et $J$ d\'efinissent des familles de sous-sch\'emas 
d'un espace affine. On calcule alors le lieu d'incidence par divisions
en utilisant le fait que les g\'en\'erateurs de $I^1$ forment une base
de Gr\"obner unitaire pour l'ordre homog\`ene avec $a>b>d>c$.
\nl
Plut\^ot que d'\'ecrire les divisions sous la forme traditionnelle 
$$f=\sum \lambda_i f_i +R$$
o\`u les $f_i$ sont les g\'en\'erateurs de $I^1$, 
nous travaillerons dans 
$$(k[[u,u',u'',v,v',v'',e,f,g,h,i,j,l,m]]\; \ox \;k[a,b,c,d])/I^1$$ et nous \'ecrirons $f=R$. Moyennant cette convention, les divisions des g\'en\'erateurs
de $J$ par $I^1$ s'\'ecrivent:
$$cv+u-a=cv+u -ec -f=(u-f)+c(v-e)$$
donc $I(W)\supset (e-v,u-f)$
$$dv+w-b=lcv+mv+w-gc-h$$
donc $I(W) \supset  (g-lv,h-w-vm)$.
$$d+cv'+u'-ac=lc+m+cv'+u'-fc-eic-ej$$
donc  $I(W) \supset (l+v'-f-ei,m+u'-ej)$.
\beq
2cd+cv''+u''-ac^2&=&2lc^2+2cm+cv''+u''-ejc+(-l-v')c^2\\
&=& (l-v')(ic+j)+2cm+cv''+u''-ejc
\eeq
donc $I(W)\supset (i(l-v')+2m+v''-ej,j(l-v')+u'')$.\\
Finalement, on a donc 
\beq
I(W)&=&(e-v,u-f,g-lv,h-w-vm,l+v'-f-ei,m+u'-ej,\\
&&i(l-v')+2m+v''-ej,j(l-v')+u'')
\eeq
et 
\beq
W&=&\s k[[u,u',u'',v,v',v'',e,f,g,h,i,j,l,m]]/I(W)\\
&=& \s k[[u',v',e,f,i,j]]
\eeq
$W$ est non singulier de dimension $3.dim(X)$.
\findem

\bprop
\label{R_{123}^{123} admissible}
Soit $X$ une vari\'et\'e lisse irr\'eductible. Le sch\'ema 
$R_{123}^{123}(X)$ est une vari\'et\'e lisse irr\'eductible
de dimension $3.dim(X)$.
En particulier $\eta_{123}^{123}$ est admissible.
\eprop
on fait ici aussi la d\'emonstration dans le cas 
o\`u $X$ est de  dimension deux.
\nl
Soit $(p^{123},p_{123})$ un point sp\'ecial de $R_{123}^{123}(X)$. 
Notons $p_{12}$ le point de $\hdx$ d\'efini par le support de $[p^{123}]$.
On reprend les notations de la d\'emonstration de la proposition 
\ref{R_{123}^1 admissible} pour les voisinages formels des points 
$p_{12},p_{123}$   ainsi que pour l'id\'eal $J$.
\nl 
Le voisinage formel de $p^{123}$ est isomorphe \`a 
$$\s k[[e,f,g,h,i,j,l,m,n,o,p,q]]$$
et l'id\'eal universel de $Hilb^{3,2}(X) \x \hdx$ au dessus de cette carte 
est 
$$I^{123}=(a-ec^2-fc-g,b-hc^2-ic-j,c^3-lc^2-mc-n,d-oc^2-pc-q)$$
Il nous faut calculer le lieu $W$ de 
$$\s k[[u,u',u'',v,v',v'',e,f,g,h,i,j,l,m,n,o,p,q]]$$
au dessus duquel $I^{123} \supset J$. 
On effectue des divisions:
$$cv+u-a=cv+u-ec^2-fc-g$$
donc $e=0,\ g=u, \ f=v\ \ \mbox{ sur } W$.
$$dv+w-b=voc^2+pcv+qv+w-hc^2-ic-j$$
donc $vo=h,\ vp=i,\ vq+w=j \ \ \mbox{ sur } W$.
$$d+cv'+u'-ac=oc^2+pc+q+v'c+u'-fc^2-gc$$
donc $o=f,\ p+v'=g,\ q+u'=0 \ \ \mbox{ sur } W$.
\beq
2cd+cv''+u''-ac^2&=& 2oc^3+2pc^2+2qc+cv''+u''-fc^3-gc^2\\
&=& (2o-f)(lc^2+mc+n)+(2p-g)c^2+c(2q+v'')+u''
\eeq
donc $$I(W)\supset ((2o-f)l+2p-g,\ (2o-f)m+2q+v'',\ (2o-f)n+u'')$$
Au total 
\beq
I(W)=(e,\ g-u, \ f-v, \ vo-h,\ vp-i,\ vq+w-j,\ o-f,\ p+v'-g,\\ 
q+u',\ (2o-f)l+2p-g,\ (2o-f)m+2q+v'',\ (2o-f)n+u'')
\eeq
et 
\beq
W&=&\s k[[u,u',u'',v,v',v'',e,f,g,h,i,j,l,m,n,o,p,q]]/I(W)\\
&=& \s k[[l,m,n,o,p,q]]
\eeq
Donc $W$ est lisse de dimension $3.dim(X)$.
\findem

\subsubsection{Quelques enrichissements non admissibles}
Dans cette section, on montre que certains enrichissements ne  sont
pas admissibles. On utilise pour cela trois raisonnements diff\'erents,
l'un donnant les \'enonc\'es \ref{1,2,123 non admissible},
\ref{non adm compliquee} et 
\ref{deux doublets de doublets et non admissible}, un autre donnant 
les \'enonc\'es \ref{1,2,3,123 non admissible} et 
\ref{autres non admissibles}, et enfin un dernier pour
l'\'enonc\'e \ref{non admissibilite en niveau au moins trois}
\nl
On d\'eduit facilement du  lemme suivant un crit\`ere de non
admissibilit\'e (corollaire \ref{critere non admissibilite}).
\bprop \label{description du point general}
Soit $\eta$ un enrichissement admissible et $\hat t$ un point de 
$\re$. Les conditions suivantes sont
\'equivalentes quand $X$ est lisse irr\'eductible:
\bit
\item
il existe un point $p$ de $F(X,3)$ tel que $\hat t=f_{\eta}(p)$
\item
le sous-sch\'ema $[p_{\eta,\eta_{123}}(\hat t)]$ de $X$ admet pour support 
trois points distincts.
\eit
\eprop \\
1 \FD 2 est \'evident.\\
2 \FD 1: par d\'efinition d'un enrichissement admissible, $\hat t$ 
est limite de points de $f_{\eta}(F(X,3))$, ie. il existe $\phi:
\s k[[t]]_t \fd F(X,3)$ tel que $f_{\eta}\circ \phi$ se prolonge 
en le point special $0$ de $\s k[[t]]$ en prenant la valeur $\hat t$. 
Pour conclure, il nous suffit de montrer que 
$\phi$ se prolonge en $0$ puisqu'alors 
$\hat t=f_{\eta}(\phi(0))$. Ce prolongement existe si 
le prolongement de $f_{\eta_{123}} \circ
\phi$ est tel que  $[f_{\eta_{123}} \circ
\phi(0)]$ a pour support trois points distincts. 
Or cette affirmation est vraie puisque 
$f_{\eta_{123}} \circ \phi= p_{\eta,\eta_{123}}
\circ f_{\eta} \circ \phi$ se prolonge en $0$ par 
$p_{\eta,\eta_{123}}(\hat t)$.
\findem

\bcor
\label{critere non admissibilite}
Si un point $\hat t$ 
de $\re$ satisfait la deuxi\`eme condition mais pas la premi\`ere,
alors $\eta$ n'est pas admissible.
\ecor
Ce crit\`ere  s'applique aux trois propositions 
suivantes. 

\begin{prop}
\label{1,2,123 non admissible}
L' enrichissement $\eta_{1,2,123}$ n'est pas admissible
\end{prop}
\demo: le point 
$(p_1,p_1,p_1 \cup p_2 \cup p_3) \in \hachun \x \hachdeux \x \hachudt$
est un point de $R_{1,2,123}(X)$ qui n'est pas dans 
$f_{\eta_{1,2,123}}(F(X,3))$
\findem

\begin{prop}
\label{non adm compliquee}
Soit $\eta$ un enrichissement de niveau deux  contenant  $\sigma ^1$ 
comme  unique structure de niveau deux, contenant  
$\sigma_1$ ou $\sigma_2$. Alors $\eta$ n'est pas admissible. 
\end{prop}
\demo: consid\'erons le point  
$p=(p_2,p_1,p_3,p^1,p_{123})$
de 
$$
\hachun \x  \hachdeux \x \hachtrois \x H^1(X) \x \hachudt
$$
avec 
$$[p_{ij}]=p_i\cup p_j\ \ [p^1]=p_{12} \cup p_{13}\ \ [p_{123}]=p_1 \cup 
p_2 \cup p_3$$
La projection de $p$ sur $H_{\eta}(X)$ est un point de $R_{\eta}(X)$ qui
n'est pas dans $f_{\eta}(F(X,3))$.
\findem

\begin{prop} \label{deux doublets de doublets et non admissible}
Soit $\eta$ un enrichissement de niveau deux contenant $\sigma^1$ et 
$\sigma^2$, et  tel que $\eta\cap \linebreak[0]
\{\sigma_{12},\sigma_{13},\sigma_{23}\}\linebreak[0]=\emptyset$. Alors $\eta$ 
est non admissible.
\end{prop}
\demo: consid\'erons le point 
$$\hat t=(p_1,p_2,p_3,p^1,p^1,p^3,p^{123}, p_{1} \cup p_2 \cup p_3) $$
de 
$$\hachun \x \hachdeux \x \hachtrois \x H^1(X) \x H^2(X)
\x H^3(X) \x H^{123}(X)\x \hachudt$$  
Le point 
$$p_{\eta_{1,2,3,123}^{1,2,3,123},\eta}(\hat t)$$  
 est un point de $R_{\eta}(X)$ qui n'est pas dans $f_{\eta}(F(X,3))$. 
\findem

\bprop
\label{1,2,3,123 non admissible}
L'enrichissement $\eta_{1,2,3,123}$ n'est pas admissible
\eprop
Soit $X$ une surface lisse. On veut montrer que les 
sch\'emas 
$$\overline{f_{\eta_{1,2,3,123}}(F(X,3))}
 \inc \hachun \x \hachdeux \x \hachtrois \x \hachudt$$
et 
$$R_{1,2,3,123}(X) \inc \hachun \x \hachdeux \x \hachtrois \x \hachudt$$
sont diff\'erents. Pour cela, on les projette sur 
$ \hachun \x \hachdeux \x \hachudt$:
$$p:\overline{f_{\eta_{1,2,3,123}}(F(X,3))} \fd \hachun \x \hachdeux  \x \hachudt$$
$$q:R_{1,2,3,123}(X) \fd \hachun \x \hachdeux \x \hachudt$$
On consid\'ere un point 
$$
\hat t=(x,x,t) \in \hachun \x \hachdeux  \x \hachudt$$
o\`u $[t]$ est un deux gros point et $x \inc [t]$. On va montrer que la fibre $p^{-1}(\hat t)$ est r\'eduite \`a un 
point et que ce n'est pas le cas  de $q^{-1}(\hat t)$.
\\
Pour la deuxi\`eme  affirmation,  la fibre $q^{-1}(\hat t)$ est le lieu 
de $\hachtrois=X$ form\'e par les $x'$ tels que:
\bit
\item 
$[t] \inc x.x.x'$
\item
$x' \inc [t]$
\eit
La premi\`ere condition est v\'erifi\'ee car on a d\'ej\`a $[t]\inc x.x$
La deuxi\`eme condition dit alors que  $q^{-1}(\hat t)=[t]$,
qui n'est pas r\'eduit \`a un point.
\nl
Notons $W$ l'image du morphisme $p$.
Pour montrer la premi\`ere affirmation, il nous suffit de construire un 
morphisme 
$$W \fd \hachtrois=X$$
qui envoie un point g\'en\'eral $(p_1,p_2,p_1 \cup p_2 \cup p_3)$ de $W$ 
sur $p_3$. A l'aide de l'inclusion
$$W \inc \hachun \x \hachdeux \x \hachudt $$
on r\'ecolte alors un morphisme 
$$\phi : W \fd  \hachun \x \hachdeux \x \hachtrois \x \hachudt$$
qui envoie un point g\'en\'eral $(p_1,p_2,p_1 \cup p_2 \cup p_3)$ sur 
$(p_1,p_2,p_3,p_1 \cup p_2 \cup p_3)$.
Donc 
$$
\phi(W) \supset f_{\eta_{1,2,3,123}}(F(X,3))
$$ 
et par suite, 
$$
\phi(W) \supset \overline{f_{\eta_{1,2,3,123}}(F(X,3))}
$$ 
Les fibres du 
morphisme 
$$\overline{f_{\eta_{1,2,3,123}}(F(X,3))} \fd W$$ 
sont incluses dans les fibres du morphisme 
$$\phi(W) \fd W$$ qui sont r\'eduites \`a un point.
\\
Il nous reste \`a construire ce morphisme $W \fd X$. 
Soit $U$ l'ouvert $f_{\eta_{1,2,123}}(F(X,3))$ de $W$.
Le morphisme 
\beq
\phi:U &\fd& W\x \hdx\\
(p_1,p_2,p_{123}) &\mapsto& (p_1,p_2,p_{123},p_1 \cup p_2)
\eeq
d\'efinit une sous-vari\'et\'e 
$$\overline{\phi(U)} \inc W\x Hilb^2(X)$$
Un point $(p_1,p_2,p_{123},p_{12})$ de $\phi(U)$ v\'erifie $[p_{12}]
 \inc [p_{123}]$. Cette condition reste vraie sur $\overline{\phi(U)}$
et on peut d\'efinir
\beq 
\psi:\overline{\phi(U)} &\fd& X\\
(p_1,p_2,p_{123},p_{12}) &\mapsto & Res([p_{12}],[p_{123}])
\eeq 
Pour montrer que $\psi$ induit un morphisme $W\fd X$, il suffit de 
v\'erifier que $\psi$ est constant sur les fibres de la projection
$$\overline{\phi(U)} \fd W$$
V\'erifions le sur un point sp\'ecial de $W$. Un tel point est de la 
forme $(p,p,t)$ o\`u $[t]$ est un deux-gros point de support $p$. 
Un point de la fibre est de la forme $(p,p,t,d)$ avec $[d]\inc [t]$ et 
donc $Res([d],[t])=p$. Le morphisme $\psi$ est constant sur les 
fibres.
\findem \nl
De la m\^eme mani\`ere, on montre:
\bprop
\label{autres non admissibles}
Les enrichissements $\eta_{1,2,123}^{123}$ et $\eta_{1,2,3,123}^{123}$
ne sont pas admissibles.
\end{prop}

\begin{lm}\label{non admissibilite en niveau au moins trois}
  Soit $\eta$ un enrichissement de $\{1,2,3\}$ de niveau $\geq 3$ 
  dont les
  structures de niveau un et deux sont incluses dans l'ensemble
  $\{\sigma_1,\sigma_2, \sigma_3,\sigma_{12},\sigma_{13},
  \sigma_{123},\sigma ^1\}$. Alors $\eta_0$ n'est pas admissible.
\end{lm}
\demo: soit $p$ un point de $X$, $t$ un point de $Hilb^3(X)$ tel que $[t]$
soit  support\'e par $p$ et isomorphe \`a $\s k[x,y]/(x^2,xy,y^2)$. Soit 
$d_1,\dots,d_6$ des points de $Hilb^2(X)$ distincts 
tels que $[d]\inc [t]$. D\'efinissons $d^{ij}\in Hilb^{2,2}(X)$ par 
$[d_{ij}]=d_i \cup d_j$. Utilisons ces donn\'ees pour construire 
un point de $R_{\eta}(X)$ qui n'est pas dans $\overline{f_{\eta}(F(X,3))}$. 
Notons $\eta_1$ (resp. $\eta_2$) 
l'ensemble des structures de niveau au plus 
deux (resp. au moins trois) de $\eta$ de sorte que
$H_{\eta}(X)=H_{\eta_1}(X) \x 
H_{\eta_2}(X)$. Si $\sigma$ est une structure de niveau $l$
de $\eta_2$, 
alors $\sigma$ est dans $\Sigma_{p2\dots 2}(E)$ o\`u $p$ vaut deux ou
trois. On peut voir $\sigma$ comme une structure de niveau $l-2$ sur 
$\Sigma_{22}(E)$. Puisque $\Sigma_{22}(E)$ est un ensemble \`a trois 
\'el\'ements, il s'identifie \`a $E$ 
 par la bijection faisant correspondre $\{\{i,j\},\{i,k\}\}$ et
$\{i\}$.
La structure $\sigma$ 
peut alors \^etre identifi\'e \`a une structure $\sigma'$ de niveau 
$l-2$ de $E$.
Notons  $\eta'_2:=\cup_{\sigma\in \eta_2}\sigma'$.
Notons $q$ la projection du point 
$(p,p,p,d_1,d_2,d^{12},t)$ de 
$$ H_{\sigma_1}(X)\x H_{\sigma_2}(X)
\x H_{\sigma_3}(X) \x H_{\sigma_{12}}(X) \x H_{\sigma_{13}}(X)\x 
H_{\sigma^1}(X)
\x H_{\sigma_{123}}(X)$$ sur $H_{\eta_1}(X)$.
Puisque $H_{\sigma'}(Hilb^{22}(X))=H_{\sigma}(X)$, le point 
$P=(q,f_{\eta'_2}(d^{12},d^{34},d^{56}))$ est dans 
$H_{\eta_1}(X)\x H_{\eta'_2}(H_{22}(X))=H_{\eta}(X)$. V\'erifions 
que $P$ est dans  $R_{\eta}(X)$ mais pas dans $\overline{f_{\eta}(F(X,3))}$.
Ce point n'est pas dans $\overline{f_{\eta}(F(X,3))}$ car si $\sigma$ est une
structure de niveau au moins trois de $\eta$, le point 
$p_{\eta,\sigma}(P)$ est d\'efini \`a l'aide d'au moins quatre
doublets parmi $d_1,\dots,d_6$. Or, s'il \'etait dans 
$\overline{f_{\eta}(F(X,3))}$, il serait d\'etermin\'e par au plus 
trois doublets. Pour voir qu'il est en revanche dans 
$R_{\eta}(X)$, il faut v\'erifier les relations d'incidence. 
Par construction de $q$, toutes les relations d'incidence liant les
structures de niveau au plus deux sont v\'erifi\'ees.
Puisque la composante de $P$ sur $H_{\eta_2}(X)$ est de la forme 
$f_{\eta'}(p)$ pour un point $p$, toutes 
les relations d'incidence concernant des structures de niveau 
au moins trois sont v\'erifi\'ees. Pour les relations d'incidence,
liant des structures de niveau au plus deux et des structures de 
niveau au moins trois, elles se v\'erifient soit trivialement,
soit en utilisant le fait que les doublets $[d_i]$ sont inclus dans 
$[t]$ (nous laissons l'\'ecriture pr\'ecise de la liste des
v\'erifications
au lecteur).   
\findem

\subsubsection{Lemmes de contamination}
Cette section contient les propositions qui permettent de produire 
des enrichissements admissibles $\eta'$ \`a partir d'enrichissements 
admissibles $\eta$. 
\nl
La premi\`ere de ces propositions a d\'ej\`a \'et\'e d\'ecrite:
c'est la proposition \ref{oubli=isomorphisme residuel}.
\begin{prop}
\label{d+d=d2}
Soit $\eta$ un enrichissement admissible contenant $\sigma_{12}$ et $\sigma_{13}$. Soit $\eta'=\eta \cup \{\sigma^1\}$. L'enrichissement 
$\eta'$ est admissible et le morphisme d'oubli $p_{\eta',\eta}:R_{\eta'}(X)
 \fd \re$ 
est un isomorphisme.
\end{prop}
\demo: d'apr\`es \ref{critere d'isomorphisme d'oubli}, 
il suffit de d\'emontrer que la fibre de 
$p_{\eta',\eta}$ au dessus d'un point 
$q=q(x,d,t)$ est incluse dans un point.
Les doublets inclus dans 
le triplet 
$[t]$ forment un sous-sch\'ema de $\hdx$  isomorphe \`a $\PP^1$. La fibre 
au dessus de $q$ est le lieu $Z$ de $Hilb^{2,2}(X)$ form\'e par les doublets 
de doublets $d^2$ v\'erifiant 
\beq
[d^2] \inc \PP^1\ \ (car \ \sigma^{1} \inc \sigma_{123})\\
\ [d^2] \inc d.d \ \   (car \ \sigma^{1} \inc \sigma_{12}\cup \sigma_{13})
\eeq
Si $Y \inj X$ est un sous-sch\'ema ponctuel de degr\'e $n$, le lieu de 
$Hilb^n(X)$ form\'e par les $p$ tels que $[p] \inc Y$ est sch\'ematiquement 
r\'eduit \`a un point. Pour montrer la proposition, il nous suffit 
donc de montrer que $\PP^1 \cap d.d$ d\'efinit un sous-sch\'ema 
de degr\'e deux. Mais c'est clair car c'est un sous-sch\'ema d'une courbe 
lisse d\'efini par le carr\'e d'un id\'eal maximal.
\findem
\\
Les propositions suivantes 
se d\'emontrent \'egalement par v\'erification au point 
sp\'ecial.

\begin{prop}
\label{d+d3=d2}
Soit $\eta$ un enrichissement admissible contenant $\sigma_{12}$ et $\sigma^{123}$. Soit $\eta'=\eta \cup \{\sigma^3\}$. L'enrichissement 
$\eta'$ est admissible et le morphisme d'oubli $p_{\eta',\eta}:R_{\eta'}(X)
 \fd \re$ 
est un isomorphisme.
\end{prop}

\begin{prop}
\label{d+d2=d3}
Soit $\eta$ un enrichissement admissible contenant $\sigma_{12}$ et $\sigma^{3}$. Soit $\eta'=\eta \cup \{\sigma^{123}\}$. L'enrichissement 
$\eta'$ est admissible et le morphisme d'oubli $p_{\eta',\eta}:R_{\eta'}(X)
 \fd \re$ 
est un isomorphisme.
\end{prop}

\noindent
Les deux propositions suivantes expliquent que certains enrichissements 
$\eta'$ sont admissibles car on peut trouver un enrichissement admissible
$\eta$ tel que $R_{\eta'}(X)$ s'identifie \`a une famille universelle 
au dessus de $\re$.

\begin{prop}
\label{univ en niveau un}
Soit $\eta=\{\theta_1,\dots ,\theta _s\}$ un enrichissement admissible 
et $$\eta'=\eta \cup \{\sigma\}$$ avec $$\sigma\in \{\sigma_1,\sigma_2,
\sigma_3\}$$
Supposons qu'il existe des relations d'incidence liant 
$\sigma$ aux $\theta_i$ parmi lesquelles  
$\sigma \inc \theta_i$ o\`u 
$\theta_i$ est une structure de niveau un.
Supposons en outre que la fibre du morphisme d'oubli
$$f_{\eta'}(F(X,3)) \fd f_{\eta}(F(X,3))$$ 
au dessus d'un point g\'en\'eral $(p_1,\dots,p_s)$ s'identifie
au sous-sch\'ema de $H_{\sigma}(X)$ form\'e par les $p_{\sigma}$ tels 
que $p_{\sigma} \inc [p_{i}]$. 
Alors $\eta'$ est admissible.
\end{prop}
\demo:
on peut d\'efinir trois ferm\'es de 
$$H_{\eta'}(X)=H_{\theta_1}(X)\x \dots\x H_{\theta_s}(X)\x H_{\sigma}(X)$$
\bit
\item 
le ferm\'e $\overline{f_{\eta'}(F(X,3))}$
\item
le ferm\'e $R_{\eta'}(X)$ form\'e par les 
$$(p_1,\dots,p_s,p_{\sigma})$$
satisfaisant  certaines relations d'incidence not\'ees
$\fxr_1,\dots,\fxr_r$
\item
le ferm\'e $Z$ qu'on d\'efinit de la fa\c con suivante. C'est le sous-sch\'ema 
contenant l'ensemble des points 
$(p_1,\dots,p_s,p_{\sigma})$
soumis \`a des  relations d'incidence. Ces relations d'incidence sont toutes 
les relations d'incidence pr\'ec\'edentes qui ne mettent pas en jeu 
le point $p_{\sigma}$ et la relation $[p_{\sigma}] \inc [p_i]$. C'est donc 
un sous-ensemble de relations de $\fxr_1,\dots,\fxr_r$.
\eit
Ces sous-sch\'emas satisfont les inclusions
$$\overline{f_{\eta'}(F(X,3))}\inc R_{\eta'}(X) \inc Z$$
Le sch\'ema $Z$ se projette sur $\re$ par un morphisme $p$. Notons 
$U$ l'ouvert dense $f_{\eta}(F(X,3))$ de $R_{\eta}(X)$.
Par hypoth\`ese on a l'\'egalit\'e
$$Z \cap p^{-1}(U) =f_{\eta'}(F(X,3))$$
Comme $Z$ est plat sur $\re$, on obtient
$$Z=\overline{Z \cap p^{-1}(U)}=\overline {f_{\eta'}(F(X,3))}$$
Par suite 
$$R_{\eta'}(X)=\overline {f_{\eta'}(F(X,3))}$$
et $\eta'$ est admissible.
\findem

\begin{prop}
\label{univ en niveau 2}
Soit $\eta=\{\theta_1,\dots ,\theta _s\}$ un enrichissement admissible 
et $$\eta'=\eta \cup \{\sigma\}$$ avec $$\sigma\in \{\sigma_{12},\sigma_{13},
\sigma_{23}\}$$
Supposons qu'il existe des relations d'incidence liant 
$\sigma$ aux $\theta_i$ parmi lesquelles  
$\sigma \inc \theta_i $ o\`u $\theta_i$ est une structure de niveau deux.
Supposons en outre que la fibre du morphisme d'oubli
$$f_{\eta'}(F(X,3)) \fd f_{\eta}(F(X,3))$$ 
au dessus d'un point g\'en\'eral $(p_1,\dots,p_s)$ s'identifie
au sous-sch\'ema de $H_{\sigma}(X)$ form\'e par les $p_{\sigma}$ tels 
que $p_{\sigma} \inc [p_{i}]$.
Alors $\eta'$ est admissible.
\end{prop}
\demo:
comme ci-dessus, $\retaprime$ s'identifie \`a une vari\'et\'e universelle 
et peut \^etre obtenue comme adh\'erence de
$ f_{\eta'} (F(X,3)) $.
\findem

\bprop \label{ajouter des 3222}
Soit $\eta$ un enrichissement admissible de $E=\{1,2,3\}$,
 $s_{l}$ l'unique structure de
$\Sigma_{3\underbrace{22\dots 
  2}_{l-1\  fois}}(E)$. 
Soit $l_0\geq 1$ le plus grand entier tel que $s_{l_0} \in \eta$.
Si $l_0 \geq 2$ alors $\eta':=\eta \cup \{s_{l_0+1}\}$  
est admissible et le morphisme d'oubli $p_{\eta',\eta}$ est un 
isomorphisme.
\eprop
pour $2\leq i \leq l_0-1$, posons 
$\eta_i:=\eta-\{s_{l_0},s_{l_0-1},\dots,s_{i+1}\}$,
$V_{l}:=p_{\eta,\eta_l}(R_{\eta}(X))$. Notons
$V_{l_0}=R_{\eta}(X)$, $V_{l_0+1}=R_{\eta'}(X)$, $p_{lk}$ le 
morphisme d'oubli de $V_l$ dans $V_k$ pour $2\leq k\leq l\leq l_0+1$. 
Soit $p_{l_0}$ un point de $V_{l_0}$, $p_s:=p_{l_0,s}(p_{l_0})$.
On suppose que $[p_{\eta,\eta_{123}}(p_{l_0})]$ est  isomorphe \`a 
$Spec\ k[x,y]/(x^2,xy,y^2)$.
On va montrer que $[p_{\eta,\{s_2\}}(p_{l_0})]$ est curviligne,
puis que si $[p_{\eta,\{s_l\}}(p_{l_0})]$ est curviligne,
alors  la fibre de $p_{l+1,l}$ au dessus de $p_l$ est 
sch\'ematiquement incluse dans    un point 
et  $[p_{\eta,\{s_{l+1}\}}(p_{l_0})]$ est aussi curviligne.
Ces r\'esultats impliqueront par r\'ecurrence que la fibre de 
$p_{l_0+1,l_0}$ au dessus de $p_{l_0}$
est incluse dans un point,  ce qui suffit pour \'etablir la proposition
d'apr\`es \ref{critere d'isomorphisme d'oubli}. 
Tout d'abord, $p_{\eta,\{s_2\}}(p_{l_0})$ est inclus dans 
$p_{\eta,\eta_{123}}(p_{l_0})$ puisque $\{s_2\}\inc \sigma_{123}$.
Puisque les doublets inclus dans $[p_{\eta,\{\eta_{123}\}}(p_{l_0})]$
forment un $\PP^1$, cela implique que
$[p_{\eta,\{s_2\}}(p_{l_0})]$
est curviligne. Un calcul en coordonn\'ees locales 
montre que si $Z$ est un sous-sch\'ema ponctuel 
de degr\'e trois curviligne d'une vari\'et\'e  $Y$, alors le lieu 
sch\'ematique de $Hilb^2(Y)$ param\'etrant les doublets $d$ 
inclus dans $Z$ est \'egalement un triplet ponctuel curviligne.
Le lieu de $Hilb^{3,2}(Y)$ param\'etrant les triplets de doublets 
inclus dans $Z$ est alors r\'eduit sch\'ematiquement \`a un point.
En particulier, si $Z=[p_{\eta,\{s_l\}}(p_{l_0})]$ est curviligne,
la fibre $F$ de $p_{l+1,l}$ vue comme sous-sch\'ema de 
$H_{s_{l+1}}$  est  sch\'ematiquement form\'ee
de triplets 
de doublets inclus dans $Z$, donc est r\'eduite \`a au plus un point.
Pour $l<l_0$, $F$ est non vide car elle
contient $p_{\eta,\{s_{l+1}\}}(p_{l_0})$
Donc  $[p_{\eta,\{s_{l+1}\}}(p_{l_0})]$ est curviligne par 
ce qui pr\'ec\`ede. 
\findem

\begin{prop} \label{ajouter des 2222}
  Soit $\eta$ un enrichissement admissible contenant deux structures 
$e\in \Sigma_{322\dots 2}$ et $g \in \Sigma_{122\dots 2}$ de m\^eme 
niveau $l$. Alors $[g]\subset [e]$. En outre, si $f$ est d\'efini
par $[f]=[e]-[g]$, alors $\eta':=\eta \cup \{f\}$ est admissible et 
le morphisme d'oubli $p_{\eta',\eta}$ est un isomorphisme. 
\end{prop}
\demo: l'inclusion $[g]\subset [e]$ se montre ais\'ement par r\'ecurrence 
sur le niveau. Le m\^eme raisonnement que dans la proposition
pr\'ec\'edente montre que la fibre sp\'eciale du morphisme d'oubli
s'identifie \`a un lieu $L\inc Hilb^2(W)$ param\'etrant 
des doublets $d$ v\'erifiant $[d]\inc [t]$ et $I([d])I([p])\inc
I([t])$, o\`u $[t]\inc W$ est un triplet curviligne ponctuel et 
o\`u $p$ est le support de $[t]$. Mais le r\'esiduel d'un sch\'ema
dans un sch\'ema curviligne \'etant bien d\'efini, $[d]$ est 
n\'ecessairement le r\'esiduel de $p$ dans $[t]$ et la fibre 
est sch\'ematiquement r\'eduite \`a un point. \findem

\subsubsection{Classification en niveau au plus deux}
Dans cette section, on classifie tous les compactifications 
de niveau au plus deux. Pour cela, on commence par produire 
une liste d'enrichissements admissibles:
\begin{prop}
\label{liste d'admissibles}
Les enrichissements
$$
\{\eta_{123},\ \eta_{1,123},\ \eta_{1,2,3,12,123},\ \eta_{3,12,123}$$
$$\eta_{123}^{1},\ \eta_{1,2,3,12,13,123}^{1},\ \eta_{123}^{123}$$
$$\eta_{1,123}^{123},\ \eta_{3,12,123}^{3,123},
\ \eta_{1,2,3,12,123}^{3,123},\ \eta_{max}\}$$
sont admissibles
\end{prop}
\noindent
On montre ensuite qu'on a ainsi toutes les compactifications 
\`a isomorphisme de compactifications pr\`es: 

\bprop \label{classification en niveau au plus deux}
Soit $\eta$ un enrichissement admissible de niveau au plus deux. 
Il existe une permutation
$g$ de $\{1,2,3\}$, il existe un enrichissement $\theta$ 
parmi les onze 
enrichissements du corollaire pr\'ec\'edent
tel que $\theta \supset g.\eta$ et  
tel que le morphisme d'oubli $p_{\theta ,g.\eta}$ soit un 
isomorphisme.
\end{prop} 
\noindent
\textit{D\'emonstration de la proposition \ref{liste d'admissibles}}:
soient $\eta \inc \eta'$  deux enrichissements avec $\eta$
admissible. 
Supposons que l'on puisse 
trouver une suite d'enrichissements
$$\eta_0:=\eta \inc \eta_1 \dots \inc \eta_p=\eta'$$ 
telle que 
$$\eta_{i+1}=\eta_i \cup \{\theta _{i+1}\}$$ 
o\`u $\eta_{i},\eta_{i+1}, \theta_{i+1}$ 
jouent le r\^ole de $\eta,\eta',\sigma$
d'une des propositions 
\ref{oubli=isomorphisme residuel},\ref{d+d=d2},\ref{d+d3=d2},
\ref{d+d2=d3},\ref{univ en niveau un},\ref{univ en niveau 2}. 
Alors $\eta'$ est admissible. 
On dira que $\eta'$ domine $\eta$ 
et que la suite $\theta_1, \dots,\theta_p$ est une suite de domination.
\bit 
\item
${\eta_{123}}$ est admissible car $R_{123}(X)=\htrx$ \'etant irr\'eductible, l'ouvert 
$f_{\sigma_{123}}(F(X,3))$ est dense.
\item
$\eta_{1,123}$ domine $\eta_{123}$ ($\sigma_1$ est une suite de domination)
\item
$\eta_{1,2,3,12,123}$ domine $\eta_{12,123}$ (suite:  $\sigma_1,\sigma_2,\sigma_3
$) qui est admissible par \ref{R_{12,123}admissible}
\item
$\eta_{3,12,123}$ domine $\eta_{12,123}$ (suite: $\sigma_3$)
qui est admissible par \ref{R_{12,123}admissible}
\item
$\eta^1_{123}$ est admissible par \ref{R_{123}^1 admissible}
\item
$\eta^1_{1,2,3,12,13,123}$ domine $\eta^1_{123}$ (suite: $\sigma_{12},\sigma_{13},
\sigma_2,\sigma_3,\sigma_{1}$)
\item
$\eta^{123}_{123}$ est admissible par \ref{R_{123}^{123} admissible}
\item
$\eta^{123}_{1,123}$ domine $\eta_{123}^{123}$ (suite: $\sigma_1$)
\item
$\eta_{3,12,123}^{3,123}$ domine $\eta_{123}^{123}$ (suite: $\sigma_{12},\sigma_3,\sigma^3$)
\item
$\eta_{1,2,3,12,123}^{3,123}$ domine $\eta_{3,12,123}^{3,123}$ (suite:
$\sigma_1,\sigma_2$)
\item
$\etamax$ domine $\eta^{123}_{123}$ (suite: $\sigma_{12},\sigma^3,
\sigma_{13},\sigma_{23},\sigma_1,\sigma_2,\sigma_3, \sigma^1,\sigma^2$)
\eit
\findem

\noindent
\textit{D\'emonstration de la propostion 
\ref{classification en niveau au plus deux}}: 
choisir un enrichissement de $\{1,2,3\}$ de niveau 
inf\'erieur ou \'egal \`a deux, c'est choisir un sous-ensemble 
de 
$$\{\sigma_1,\sigma_2,\sigma_3,\sigma_{12},\sigma_{13}, \sigma_{23}, 
\sigma^1,\sigma^2, \sigma^3,\sigma^{123},\sigma_{123}\}$$
contenant $\sigma_{123}$. Il y a donc $2^{10}$ tels enrichissements
et il nous faut d'abord les classifier pour pouvoir 
montrer la proposition. Nous 
allons le faire suivant les doublets d\'etermin\'es par l'enrichissement.
\nl
Commen\c cons par la remarque suivante. Soit $E$ un ensemble fini et 
$F$ un sous-ensemble de l'ensemble $\pdee$ des parties de $E$.
Le sous-ensemble $G:=\cup_{f\in F}f$ de $E$ est naturellement
stratifi\'e par la partition  la plus grossi\`ere
$P\inc \pdee$ pour laquelle:
$$\forall f \in F,\ f=\cup_{p_i \in P}\ p_i.$$

Consid\'erons maintenant $E=\{\sigma_{12},\sigma_{13},\sigma_{23}\}$
et
\beq
\pdee =  \{\sigma_{12},\sigma_{13},\sigma_{23},\sigma^{1},\sigma^2,\sigma^3,
\sigma^{123}\}.
\eeq
Un enrichissement $\eta$ d\'etermine un ensemble $F_{\eta}$ de $\pdee$: 
$$F_{\eta}:=\eta \cap \pdee\ \inc \pdee$$
donc un ensemble partitionn\'e $G_{\eta}$ inclus dans $E$.
Sept cas sont possibles:
\bit 
\item
cas un: $G_{\eta}=\emptyset$ ($\eta$ ne d\'etermine aucun doublet)
\item
cas deux: $G_{\eta}$ est un singleton ($\eta$ d\'etermine un doublet)
\item
cas trois: $G_{\eta}$ contient deux \'el\'ements et 
la partition est grossi\`ere 
($\eta$  d\'etermine deux doublets indistinguables)
\item
cas quatre: $G_{\eta}$ contient deux \'el\'ements et la partition divise 
$G_{\eta}$ en deux singletons ($\eta$  d\'etermine deux doublets ordonn\'es)
\item
cas cinq: $G_{\eta}$  contient trois \'el\'ements et la partition est 
grossi\`ere
\item
cas six: $G_{\eta}$ contient trois \'el\'ements et la partition divise 
$G_{\eta}$ en une paire et un singleton
\item
cas sept: $G_{\eta}$ contient trois \'el\'ements et la partition divise 
$G_{\eta}$ en trois singletons
\eit
Disons que deux enrichissements $\eta\inc \eta'$ sont identifiables 
si le morphisme $p_{\eta',\eta}$ est un isomorphisme. Nous allons montrer
en examinant les $\eta$ cas par cas que, si $\eta$ est admissible de
niveau au plus deux, on peut ajouter des structures \`a $\eta$ et 
obtenir un des  enrichissements mod\`eles $\eta'$ de la proposition  
\ref{liste d'admissibles}, identifiable \`a $\eta$.
Dans ce qui suit, les couples
d'enrichissements identifiables sont donn\'es par les propositions 
\ref{oubli=isomorphisme residuel}, \ref{d+d=d2}, \ref{d+d3=d2},
\ref{d+d2=d3}.\\
Remarquons qu'un enrichissement $\eta$ s'\'ecrit
$$\eta=F_{\eta}\cup \{\sigma_{123}\} \cup E_{\eta}$$
o\`u $E_{\eta}$ est un sous-ensemble de
$\{\sigma_1,\sigma_2,\sigma_3\}$. 
\nl
\textit{Cas un}: 
$\feta$ est vide donc, \`a permutation pr\`es des indices, $\eta$ est l'un 
des enrichissements suivants:
$$\eta_{123},\eta_{1,123},\eta_{1,2,123},\eta_{1,2,3,123}$$
Les deux premiers enrichissements font partie des enrichissements 
mod\`eles.
Le troisi\`eme enrichissement n'est pas admissible par 
\ref{1,2,123 non admissible}. Le dernier n'est pas non plus admissible 
par \ref{1,2,3,123 non admissible}.

\nl
\textit{Cas deux}: on a \`a l'action du groupe pr\`es 
$F_{\eta}=\{\sigma_{12}\}$
et on peut 
supposer que $\eta \supset \sigma_3$ car $\eta$ et $\eta \cup \sigma_3$ 
sont identifiables. On a alors deux cas.
Si $\eta$ ne contient ni $\sigma_1$, ni $\sigma_2$, alors 
$\eta=\eta_{3,12,123}$ qui est un des enrichissements mod\`eles.
Si $\eta$ contient $\sigma_1$ ou $\sigma_2$, par exemple $\sigma_1$ par 
sym\'etrie. Alors $\eta$ est identifiable \`a $\eta \cup \sigma_2$
$=$ $ \eta_{1,2,3,12,123}$ qui est l'un des enrichissements mod\`eles.
\nl
\textit{Cas trois}: 
on a \`a permutation pr\`es $F_{\eta}=\sigma^1$. 
Soit $\eta=\eta_{123}^1$ qui est l'un des enrichissements mod\`eles. Soit 
$\eta$ contient l'une des structures $\sigma_1,\sigma_2,\sigma_3$ et 
dans ce cas, $\eta $ n'est pas admissible. En effet, on peut 
supposer \`a sym\'etrie pr\`es que 
$\eta \supset \sigma_1$ ou $\eta \supset \sigma_2$ et $\eta$ n'est
alors pas admissible par \ref{non adm compliquee}. 
\nl
\textit{Cas quatre}: 
\`a permutation pr\`es, les seules possibilit\'es pour $F_{\eta}$ sont 
$$(F_{\eta}\supset \{\sigma_{12},\sigma^1\}\ \ ou\ \ F_{\eta}\supset \{\sigma_{12},
\sigma_{13}\})$$
et 
$$F_{\eta} \inc \{\sigma_{12},\sigma_{13},\sigma^1\}$$
Mais quitte \`a remplacer $\eta$ par un enrichissement qui lui est 
identifiable, on peut supposer 
dans le premier cas que $\eta \ni \sigma_{13}$ et dans le deuxi\`eme 
cas que $\eta \ni \sigma^{1}$. On est donc ramen\'e \`a 
$F_{\eta}=\{\sigma_{12} ,
\sigma_{13}, \sigma^1\}$. Dans ce cas, $\eta$ est identifiable \`a 
$\eta^1_{1,2,3,12,13,123}$ qui est l'un des enrichissements mod\`eles.
\nl
\textit{Cas cinq}: 
on a $F_{\eta}=\{\sigma^{123}\}$ et, \`a sym\'etrie pr\`es, 
$\eta$ est l'un des enrichissements suivants:
$$
\eta_{123}^{123},\eta_{1,123}^{123},\eta_{1,2,123}^{123},\eta_{1,2,3,123}^{123}
$$
Les enrichissements 
$\eta_{123}^{123}$ et $\eta_{1,123}^{123}$ font partie des enrichissements 
mod\`eles tandis que 
$\eta_{1,2,123}^{123}$ et $\eta_{1,2,3,123}^{123}$  ne sont 
 pas admissibles par \ref{autres non admissibles}.
\nl
\textit{Cas six}:
on peut supposer que le doublet d\'etermin\'e par la partition est 
$\sigma_{12}$ et que le doublet de doublets est $\sigma^3$. On a alors
$F_{\eta}\inc \{\sigma_{12},\sigma^3,\sigma^{123}\}$ et en outre l'une 
des trois situations suivantes:
\bit
\item
$\eta \supset \{\sigma_{12},\sigma^3\}$
\item
$\eta \supset \{\sigma_{12},\sigma^{123}\}$
\item
$\eta \supset \{\sigma^{123},\sigma^3\}$
\eit
Dans chacune des situations, on peut trouver un enrichissement $\eta'$
identifiable \`a $\eta$ avec $F_{\eta'}=\{\sigma_{12},\sigma^{123},\sigma^3\}$
et $\eta'\ni \sigma_3$. \\
Si $\eta'$ ne contient ni $\sigma_1$, ni $\sigma_2$ alors $\eta'$ est  
l'enrichissement mod\`ele $\eta_{3,12,123}^{3,123}$. Dans le cas contraire,
on peut supposer par sym\'etrie que $\eta'$ contient $\sigma_1$ et $\eta'$ est 
alors identifiable \`a l'enrichissement mod\`ele $\eta_{1,2,3,12,123}^{3,123}$.
\nl
\textit{Cas sept}:
\bit
\item
Si $\eta$ ne contient aucun doublet de doublets, ie. si $F_{\eta}
\cap \{\sigma^1,\sigma^2,\sigma^3\}=\emptyset$, alors  $\eta $ contient  
au moins deux  doublets, par exemple  $\sigma_{12},\sigma_{13}$. Il est donc 
identifiable
\`a un enrichissement contenant $\sigma^1(*)$. On peut donc supposer que 
$\eta $ contient $\sigma^1$. Maintenant, soit $\eta \supset \sigma_{23}$, soit
$\eta \supset \sigma^{123}$. Mais dans le deuxi\`eme cas, on peut
faire une nouvelle identification pour avoir $\eta\supset
\sigma_{23}$. 
Donc $\eta $ est identifiable \`a un enrichissement 
contenant $\sigma_{12},\sigma_{13}$ et $\sigma_{23}$. Mais tout 
enrichissement admissible 
contenant les trois doublets est identifiable \`a $\etamax$.
\item
Si $\eta$  contient exactement un doublet de doublets, par 
exemple  si $F_{\eta}=\{\sigma^1\}$. A sym\'etrie pr\`es, on peut supposer 
que $\eta$ contient $\sigma_{12}$ et alors $\sigma_{13}$ est un r\'esiduel. On s'est 
donc ramen\'e \`a la situation $(*)$ de ci-dessus.
\item
Si $\eta$ contient au moins deux doublets de doublets, par exemple 
$\sigma^1$ et $\sigma^2$. Alors $\eta $ contient au moins une 
structure $\sigma$ parmi $\{\sigma_{12},\sigma_{13},\sigma_{23}\}$
d'apr\`es \ref{deux doublets de doublets et non admissible}.
Par r\'esiduel, on 
retrouve alors les deux autres doublets \`a partir de 
$\sigma^1$, $\sigma^2$ et $\sigma^3$. Donc  on peut supposer que 
$\eta$ contient les 
trois doublets. Il est  alors identifiable 
\`a $\etamax$.  
\eit
\findem

\subsubsection{Classification en niveau sup\'erieur \`a deux}
On a classifi\'e dans la section pr\'ec\'edente les compactifications
de $F(X,3)$ et de ses quotients 
de niveau au plus deux. On montre dans cette section qu'on avait en fait
d\'ej\`a toutes les compactifications, ind\'ependamment du niveau. 
L'\'enonc\'e pr\'ecis suivant
ach\`eve de d\'emontrer le th\'eor\`eme \ref{thm: classification}. 

\begin{prop}\label{classification en niveau superieur}
  Pour tout enrichissement $\eta$ admissible de niveau au moins trois
  , il existe deux
  enrichissements admissibles 
  $\eta_1$ et $\eta_2$ o\`u $\eta_1$ est de niveau au plus deux  tels
  que  $\eta
  \subset \eta_2$, $\eta_1\subset \eta_2$, 
et tels que les morphismes d'oubli
  $p_{\eta_2,\eta}$ et $p_{\eta_2,\eta_1}$ soient des isomorphismes de
  compactificaction. 
\end{prop}
\demo: la d\'emarche consiste \`a prendre un 
enrichissement $\eta_2$ ``maximal'' en un sens \`a pr\'eciser. 
On montre que ce $\eta_2$ contient n\'ec\'essairement 
certaines structures. Puis on supprime dans un ordre bien choisi
les structures de $\eta_2$ pour obtenir $\eta_1$ de sorte que les 
morphismes d'oubli successifs soient des isomorphismes. 
 Plus pr\'ecis\'ement,  
soit $\eta_2$ un enrichissement de m\^eme 
niveau $l_{\eta}$ que $\eta$, contenant $\eta$, tel que $p_{\eta_2,\eta}$
soit un isomorphisme et qui soit maximal (au sens de l'inclusion)
pour ces propri\'et\'es. Si l'ensemble des structures de niveau 
au plus deux de $\eta_2$ est (\`a permutation de l'ensemble
$\{1,2,3\}$ pr\`es) inclus dans
$\{\sigma_1,\sigma_2,\sigma_3,\sigma_{12},\sigma_{13},\sigma_{123},
\sigma ^1\}$, alors $\eta$ n'est pas admissible d'apr\`es 
\ref{non admissibilite en niveau au moins trois}. 
Donc, d'apr\`es la d\'emonstration de 
\ref{classification en niveau au plus deux}, on en d\'eduit 
que $\eta_2$ d\'etermine trois doublets et que $\sigma ^{123}
\in \eta_2$. Puisque $E=\{1,2,3\}$ est de cardinal trois, une r\'ecurrence
\'evidente montre que $\Sigma_{p22\dots 2}(E)$ est de cardinal $3$ si 
$p=2$, $1$ si $p=3$ et $0$ si $p>3$. 
De plus, de la proposition \ref{ajouter des 3222}, 
on d\'eduit que $\eta_2$
contient l'\'el\'ement $e$ de $\Sigma_{322\dots2}(E)$
de niveau $l_{\eta}$. 
Si $\eta_2$ contient un \'el\'ement $f$ de $\Sigma_{222\dots 2}(E)$
de niveau $l_{\eta}$, alors $\eta_2$ contient $g 
\in \Sigma_{122\dots 2}(E)$ o\`u $[g]\cup [f]=[e]$ d'apr\`es
\ref{oubli=isomorphisme residuel}. On en d\'eduit par \ref{d+d3=d2}
que $p_{\eta,\eta-\{f\}}$ est un isomorphisme. En r\'ep\'etant
l'op\'eration avec tous les \'el\'ements $f_i$ de $\Sigma_{222\dots 2}(E)$
de niveau $l_{\eta}$, on obtient un isomorphisme
$p_{\eta,\eta-\{f_1,\dots,f_s\}}$ o\`u $e$ est l'unique structure de
niveau $l_{\eta}$ de 
$\eta-\{f_1,\dots,f_s\}$. L'application de \ref{ajouter des 3222}
montre ensuite que  
que si on note $\eta'_2=\eta-\{f_1,\dots,f_s,e\}$ 
l'ensemble des structures de niveau $l_{\eta}-1$
de $\eta_2$, alors le morphisme d'oubli $p_{\eta_2,\eta'_2}$ est un
isomorphisme. En r\'ep\'etant cette proc\'edure, on montre que si 
$\eta''_2$ est l'ensemble des structures de niveau $l_{\eta}-2$
de $\eta_2$, alors $p_{\eta,\eta''_2}$ est un isomorphisme. 
De proche en proche, on obtient 
finalement que si $\eta_1$ est l'ensemble
des structures de niveau au plus deux de $\eta_2$, alors
$p_{\eta_2,\eta_1}$ est un isomorphisme. 
\findem 

\subsection{D\'emonstration du th\'eor\`eme \ref{thm: stratification
    et oubli compatibles}}
Si $\eta$ et $\eta'$ sont deux enrichissements tels qu'il existe 
$g \in S_3$ avec $g.\eta\inc {\eta'}$, il existe un morphisme 
d'oubli $R_{\eta'}(X) \fd R_{g.\eta}(X)=R_{\eta}(X)$. 
Tous les 
morphismes du th\'eor\`eme \ref{thm: stratification
    et oubli compatibles} sont obtenus ainsi.
\\
D\'efinissons maintenant la stratification sur $R_{\eta}(X)$.
Notons $\eta(p_l,\dots,p_1):=\eta\; \cap\; \Sigma_{p_l,
\linebreak[0] \dots,p_1}(E)$,
${\cal P}_{\eta}(p_l,\dots,p_1)$ l'ensemble des parties de  
$\eta(p_l,\dots,p_1)$,
$E_{3,\eta}:=\{\sigma\in \eta \ t.q. \exists l\geq 0, \sigma\in 
\Sigma_{3\underbrace{22\dots2}_{l\ fois}}(E)\}$,
$E_{2,\eta}:=\{\sigma\in \eta \ t.q. \exists l>0, \sigma\in 
\Sigma_{\underbrace{22\dots2}_{l\ fois}}(E)$. Les strates de la
stratification sur $R_{\eta}(X)$ sont param\'etr\'ees par l'ensemble
$$Conf(\eta):=\{3,2,c,g\}^{E_{3,\eta}}\x \{2,1\}^{E_{2,\eta}}\x 
\Pi_{p_l,\dots,p_1}{\cal P}_{\eta} (p_l, \dots, p_1)$$ o\`u
$\{3,2,c,g\}^{E_{3,\eta}}$ repr\'esente les fonctions de
$E_{3,\eta}$ dans $\{3,2,c,g\}$.
Si $c=(f,g,\Pi P(p_l,\dots,p_1))\in Conf(\eta)$, et si
$\eta=\{\sigma_1,\dots,\sigma_r\}$, la strate de $R_{\eta}(X)$ correspondant 
\`a $c$ est
$S_c:=\{(p_1,\dots,p_r)$ \ t.q.:
  \begin{itemize}
  \item $\forall \sigma_i \in E_{3,\eta},\ [p_i]$ {est constitu\'e de
      trois points distincts  si } $f(\sigma_i)=3$, { de deux
      points distincts si }$f(\sigma_i)=2$, est un sch\'ema
      ponctuel curviligne si $f(\sigma_i)=c$,  un deux gros point
      si $ f(\sigma_i)=g$
  \item  $\forall \sigma_i \in E_{2,\eta},\ [p_i]$ est constitu\'e de
      deux points distincts  si $ g(\sigma_i)=2$, est un sch\'ema
      ponctuel si $g(\sigma_i)=1$
  \item $\forall \{\sigma_i,\sigma_j\}\subset
    \eta(p_l,\dots,p_1),\ p_i=p_j$ \ ssi  $\{\sigma_i,\sigma_j\}
    \subset P(p_l,\dots,p_1)$\}
  \end{itemize}
Si $\eta\inc \eta'$, l'image inverse par un morphisme d'oubli d'une
strate de $R_{\eta}(X)$ est une r\'eunion de strates de
$R_{\eta'}(X)$. En effet, on a par restriction 
une application de $\{3,2,c,g\}^{E_{3,\eta'}}\x \{2,1\}^{E_{2,\eta'}}$ dans
$\{3,2,c,g\}^{E_{3,\eta}}\x \{2,1\}^{E_{3,\eta}}$. L'intersection
avec $\eta$ donne une application ${\cal P}_{\eta'}(p_l,\dots,p_1)$
$\fd$  ${\cal P}_{\eta}(p_l,\dots,p_1)$. On peut utiliser ces
applications pour d\'efinir une 
application produit $q$ de $Conf(\eta')$ dans $Conf(\eta)$. L'\'egalit\'e 
$p_{\eta',\eta}^{-1}(S_c)=\cup_{q(c')=c}S_{c'}$ montre que l'image 
inverse d'une strate est bien une r\'eunion de strates. 
\nl
L'une des strates de $R_{\eta}(X)$ est $f_{\eta}(F(X,3))$
(correspondant \`a $f$ constante \'egale \`a $3$, \`a $g$ constante 
\'egale \`a $2$, et aux  parties $P(p_l,\dots,p_1)$ vides).
C'est une strate ouverte de fa\c con \'evidente et dense par d\'efinition 
d'un enrichissement admissible. 
\nl
Il nous reste \`a mettre en \'evidence une strate ferm\'ee incluse 
dans l'adh\'erence de toutes les autres strates. 
Pour la vari\'et\'e $R_{1,2,3,12,13,12,123}$, la strate correspondant 
\`a  $f$ constante \'egale \`a $g$, \`a  $g$ constante 
\'egale \`a $1$,  et aux  
parties $P(p_l,\dots,p_1)=\eta(p_l,\dots,p_1)$  
est form\'ee des points de la forme $q(p,d,t)$ au sens de la
proposition \ref{point special}. 
Il r\'esulte 
de [LB] que cette strate  est ferm\'ee incluse dans 
l'adh\'erence de toutes les autres strates. 
L'image inverse $S$ de cette strate sp\'eciale 
par l'isomorphisme $R_{1,2,3,12,13,12,123}(X)
\simeq R_{max}(X)$ 
est encore ferm\'ee et incluse dans l'adh\'erence de toutes 
les autres strates. Si $\eta\inc \etamax$ est un enrichissement 
quelconque, $p_{\eta_{max},\eta}(S)$ est la strate de $R_{\eta}(X)$
form\'ee des points de la forme $q(p,d,t)$. 
Elle est ferm\'ee comme image d'une vari\'et\'e projective,
incluse dans l'adh\'erence de toutes les autres strates
par continuit\'e 
que $p_{\eta_{max},\eta}(S)$. 
\findem 

\begin{rem}
  On peut en fait montrer que l'adh\'erence d'une strate de la 
stratification pr\'ec\'edente est une r\'eunion de strates. 
\end{rem}

\section{ Etude des quotients par les actions naturelles}
\label{sec:etude_des_quotients}
On a vu \`a la section \ref{sec:prop resultant des def} 
que chaque sch\'ema $\re$ est 
naturellement muni d'une action de groupe. Dans cette section nous 
d\'eterminons les quotients des sch\'emas $\re$ par ces actions
quand $\eta$ est l'un des enrichissements du
th\'eor\`eme de classification, i.e. nous montrons les 
th\'eor\`emes  \ref{thm: quotients} et \ref{thm: quotient=oubli}
de l'introduction.
\nl
L'id\'ee pour d\'eterminer les quotients $\re/G$ et  
les morphismes quotients consiste \`a  trouver pour chaque $\eta$
un enrichissement 
$\eta'\inc \eta$ tel que $R_{\eta'}(X) \simeq R_{\eta}(X)$ et tel que 
l'action du groupe sur  $R_{\eta'}(X) $
soit plus facilement ma\^ \i trisable. 
Toutes les d\'emonstrations \'etant identiques, nous ne traiterons 
qu'un seul cas.
Montrons par exemple $R^1_{1,2,3,12,13,123}(X)/S_2=R^{1}_{123}(X)$.
\\
Nous avons vu lors de la classification que $R^1_{1,2,3,12,13,123}(X)$ \'etait 
isomorphe \`a 
$$R^1_{12,13,123}(X) \inc \hachundeux \x \hachuntrois \x H^1(X) \x \hachudt$$
et l'action de $S_2$ sur cette vari\'et\'e est:
$$\varepsilon.(p_{12},p_{13},p^1,p_{123})=(p_{13},p_{12},p^1,p_{123})$$
o\`u $\varepsilon$ est l'\'el\'ement de $S_2$ diff\'erent de l'identit\'e.
\\
Quand la caract\'eristique du corps est premi\`ere au cardinal d'un 
groupe $G$, quand ce groupe $G$ agit sur un sch\'ema projectif $Y$ en laissant 
un sous-sch\'ema $Z$ invariant, le sch\'ema quotient $Z/G$ est naturellement 
un sous-sch\'ema de $Y/G$. Dans notre cas, cela veut dire qu'on a:
\begin{eqnarray*}
  R^1_{12,13,123}(X) /S_2 \inj (\hachundeux \x \hachuntrois \x & H^1(X)&
\x H_{123}(X))/S_2\\
&=&\\
 (\hachundeux \x \hachuntrois)/S_2 \x & H^1(X)&
\x H_{123}(X)
\end{eqnarray*}
car $S_2$ n'agit pas sur les deux derniers facteurs. 
\\
Il existe un morphisme naturel 
$$\phi:H^1(X) =Hilb^{2,2}(X) \fd (\hdx\x \hdx)/S_2 =(\hachundeux \x \hachuntrois)/S_2$$
d\'efinissant une sous-vari\'et\'e
$$V\inj (\hachundeux \x \hachuntrois)/S_2 \x H^1(X) \x \hachudt$$
form\'ee par les $(p,q,r)$ tels que $p=\phi(q)$. La projection 
$$V \fd H^1(X) \x H_{123}(X)$$ est un isomorphisme. 
Le quotient $R^1_{12,13,123}(X)/S_2$ est une sous-vari\'et\'e de $V$ 
comme on le v\'erifie en un point g\'en\'eral. On peut donc le projeter
par un isomorphisme sur $H^1(X) \x \hachudt$. 
En r\'esum\'e le morphisme quotient 
$$R^1_{1,2,3,12,13,123}(X) \fd H^1(X) \x \hachudt$$
s'\'ecrit comme compos\'ee:
\beq
R^1_{1,2,3,12,13,123}(X) \fd R^1_{12,13,123}(X)\inc \hachundeux \x \hachuntrois &\x& H^1(X) \x \hachudt\\
&\fb&\\
(\hachundeux \x \hachuntrois)/S_2 &\x& H^1(X)\x H_{123}(X)\\
&\fb&\\
H^1(X) &\x& \hachudt
\eeq
C'est donc simplement un morphisme d'oubli et l'image de $R^1_{1,2,3,12,13,123}(X)$
par ce morphisme d'oubli est $R^1_{123}(X)$.
\findem

\section{Comparaison aux  constructions classiques}
\label{sec identification des classiques}
Nous allons dans cette section montrer que nos 
compactifications englobent les compactifications 
de $F(X,n)$  
de Schubert-Semple-Le Barz,de  Kleiman pour $n\leq 3$, 
et de Cheah, 
ainsi que leurs quotients.
Pour cela nous nous contentons de rassembler ici
le travail d\'eja effectu\'e par Le Barz et par Keel.
L'ensemble des r\'esultats de cette section est r\'esum\'e par 
le th\'eor\`eme \ref{comparaison aux classiques}.
\\
On peut montrer en revanche
([Ev1]) 
que la compactification  $X[3]$ de Fulton-MacPherson
n'est isomorphe \`a
aucune de nos vari\'et\'es de triplets.\\
Comme il est \'evident  d'apr\`es les  d\'efinitions
que notre construction englobe les sch\'emas
\'etudi\'es par Cheah
([Ch]), nous n'aborderons ici que les 
deux autres constructions.

\subsection{La vari\'et\'e de Schubert-Semple }
Cette compactification $W_6^{*}$ de $F(\plp,3)$ 
est d\'efinie comme  la 
sous-vari\'et\'e  de 
$$\plp\x \plp\x \plp\x\plp\:^{*}\x \plp\:^{*}\x \plp\:^{*}\x \GG(2,\PP^5)$$
form\'ee par les 
$(p_1,p_2,p_3,D_{12},D_{13},D_{23},\Sigma)$
tels que 
\bit
\item
la droite $D_{ij}$ contient les points $p_i$ et $p_j$
\item
$\Sigma$ est un r\'eseau de coniques de $\plp$ contenant les coniques 
d'id\'eal $I(D_{ij}).I(D_{jk})$.
\eit
Schubert voulait \`a des fins \'enum\'eratives
une vari\'et\'e universelle lisse  param\'etrant les triangles du
plan, \'eventuellement d\'eg\'en\'er\'es.Mais une construction 
de cette vari\'et\'e dans:
$\plp\x \plp\x \plp\x\plp\:^{*}\x \plp\:^{*}\x \plp\:^{*}$ 
donnerait une vari\'et\'e 
singuli\`ere. L'introduction du r\'eseau de coniques est un moyen 
g\'eom\'etrique de d\'esingulariser cette vari\'et\'e.
\nl
Comme l'a montre Le Barz,
 $R_{1,2,3,12,13,23,123}(X)$ co\"\i ncide avec la construction
de Schubert-Semple dans le cas o\`u $X=\plp$. Le  point
$(p_1,p_2,p_3,p_{12},p_{13},p_{23},p_{123})$
de 
$X\x X\x X \x \hdx \x \hdx \x \hdx \x\htrx$
correspond au point de la construction de Schubert-Semple dont 
les  droites $D_{ij}$ sont les droites qui contiennent le
sch\'ema $[p_{ij}]$, et dont le r\'eseau de coniques est form\'e par
les coniques  
contenant $[p_{123}]$.
\nl
D'apr\`es l'\'etude des passages au quotient, le quotient de 
$R_{1,2,3,12,13,23,123}(X)$ $\simeq R_{max}(X)$
par $S_3$ est $R_{123}^{123}(X)$ et le quotient par $S_2$ est 
$R_{3,12,123}^{3,123}(X)$.

\subsection{Les vari\'et\'es de Kleiman}
Ces compactifications $X_r$ de $F(X,r)$ ont \'et\'e introduites dans [K] et elles
ont \'et\'e 
abondamment \'etudi\'ees depuis ([R],[DO],[W],[Hb] par exemple).
\nl
Leur d\'efinition est donn\'ee par r\'ecurrence sur $r$. On pose 
$X_0:=\s k$, $X_1:=X$ et 
$X_{r+1}=Bl_{\Delta} \;(X_r\x _{X_{r-1}}X_r)$, o\`u $\Delta$ 
est la diagonale. G\'eom\'etriquement, un point de $X_r$ 
correspond \`a un
$r$-uplet ordonn\'e $(p_1,\dots,p_r)$. 
Chaque  $p_i$ est soit un point de $X$, 
soit un point infiniment voisin 
( ie. un point situ\'e 
sur une surface $\tilde{X}$ obtenue 
par \'eclatement successifs de points).
\nl
Pour $r=2$, les vari\'et\'es $X_2=Bl_{\Delta}(X\x X)$ 
et $R_{1,12}(X)$ sont isomorphes. 
\nl
Keel a montr\'e l'isomorphisme  $R_{1,2,3,12,13,123}(X) \simeq X_3$
[Kee]. D'apr\`es la classification et l'\'etude des quotients, 
$X_3/S_2\simeq R_{1,2,3,12,13,123}^1(X)/S_2=R_{123}^1(X)$.

\newpage
\pagestyle{plain}
\noindent
{\bf \huge Bibliographie}
\\[10mm]
[B]: Brian\c con J., Description de $Hilb^n\CC\{x,y\}$, Invent .Math. 41,
(1977)
\nl
[CF]: Collino A., Fulton W., Intersection rings of spaces of
triangles, Colloque en l'honneur de Pierre Samuel,
M\'em. Soc. Math. France, no. 38, (1989), 75-117.
\nl
[Ch]: Cheah J.,
Cellular decompositions for nested Hilbert schemes of points,
Pacific J. Math. 183 (1998), no. 1, 39--90. 
\nl
[DO]: Dolgachev I., Ortland D., Points sets in projective spaces and 
theta functions, Ast\'erisque, 165, (1988)
\nl
[E]: Eisenbud D., Commutative Algebra with a view towards algebraic geometry,
Springer
\nl
[Ev]: Evain L., Collisions de  gros points sur une surface, 
en pr\'eparation.
\nl
[Ev1]: Evain L., Collisions de trois gros points sur une surface
algebrique, Th\`ese, Nice, (1997)
\nl
[FMP]: Fulton W., MacPherson R., 
A compactification of configuration of spaces, 
Annals of Math., 139 (1994)
\nl 
[Hb]: Harbourne B., Iterated blow-ups and moduli for rational surfaces, 
Lecture Notes, 1311, (1988)
\nl
[K]: Kleiman S., Multiple point formulas I: Iteration, Acta. Math., 147,(1981)
\nl
[Kee]: Keel S., Functorial construction of Le Barz's triangle space 
with applications, Transactions of the A.M.S., 335, (1993)
\nl 
[LB]: Le Barz P., La vari\'et\'e des triplets complets, Duke Math., 57, (1988) 
\nl
[R]: Ran Z., Curvilinear enumerative geometry, Acta. Math., 155, (1985)
\nl
[RS]: Roberts J., Speiser R., Enumerative geometry of triangles,
Comm. Algebra 12 ,(1984), no. 9-10, 1213-1255
\nl 
[Sch]: Schubert H., Anzahlgeometrische Behandlung des Dreieckes, Math. Ann. 17 (1880)
\nl 
[S]: Semple J.G., The triangle as a geometric variable, Mathematika 1,
(1954),
80-88
\nl
[SK]: Semple J.G., Kneebone G., Algebraic curves, Oxford University press,
London, (1959)
\nl 
[Ty]: Tyrell J.E., , Mathematika 6, (1959), 158-164
\nl
[W]: Walter C., Collisions of clusters of infinitely near points, non
publi\'e.

\end{document}

%% file: mor_oublis.pstex_t
\begin{picture}(0,0)%
\includegraphics{mor_oublis.pstex}%
\end{picture}%
\setlength{\unitlength}{2072sp}%
\begingroup\makeatletter\ifx\SetFigFont\undefined%
\gdef\SetFigFont#1#2#3#4#5{%
  \reset@font\fontsize{#1}{#2pt}%
  \fontfamily{#3}\fontseries{#4}\fontshape{#5}%
  \selectfont}%
\fi\endgroup%
\begin{picture}(9720,3930)(91,-4291)
\put(9811,-1456){\makebox(0,0)[lb]{\smash{\SetFigFont{6}{7.2}{\rmdefault}{\mddefault}{\updefault}\special{ps: gsave 0 0 0 setrgbcolor}$R_{123}(X)$\special{ps: grestore}}}}
\put( 91,-2311){\makebox(0,0)[lb]{\smash{\SetFigFont{6}{7.2}{\rmdefault}{\mddefault}{\updefault}\special{ps: gsave 0 0 0 setrgbcolor}$R_{max}(X)$\special{ps: grestore}}}}
\put(1216,-1636){\makebox(0,0)[lb]{\smash{\SetFigFont{6}{7.2}{\rmdefault}{\mddefault}{\updefault}\special{ps: gsave 0 0 0 setrgbcolor}$R_{1,2,3,12,13,123}^{1}(X)$\special{ps: grestore}}}}
\put(1171,-3211){\makebox(0,0)[lb]{\smash{\SetFigFont{6}{7.2}{\rmdefault}{\mddefault}{\updefault}\special{ps: gsave 0 0 0 setrgbcolor}$R_{1,2,3,12,123}^{3,123}(X)$\special{ps: grestore}}}}
\put(3826,-556){\makebox(0,0)[lb]{\smash{\SetFigFont{6}{7.2}{\rmdefault}{\mddefault}{\updefault}\special{ps: gsave 0 0 0 setrgbcolor}$R_{1,2,3,12,123}(X)$\special{ps: grestore}}}}
\put(3601,-4291){\makebox(0,0)[lb]{\smash{\SetFigFont{6}{7.2}{\rmdefault}{\mddefault}{\updefault}\special{ps: gsave 0 0 0 setrgbcolor}$R_{3,12,123}^{3,123}(X)$\special{ps: grestore}}}}
\put(6031,-1231){\makebox(0,0)[lb]{\smash{\SetFigFont{6}{7.2}{\rmdefault}{\mddefault}{\updefault}\special{ps: gsave 0 0 0 setrgbcolor}$R_{3,12,123}(X)$\special{ps: grestore}}}}
\put(5896,-2986){\makebox(0,0)[lb]{\smash{\SetFigFont{6}{7.2}{\rmdefault}{\mddefault}{\updefault}\special{ps: gsave 0 0 0 setrgbcolor}$R_{1,123}^{123}(X)$\special{ps: grestore}}}}
\put(8191,-556){\makebox(0,0)[lb]{\smash{\SetFigFont{6}{7.2}{\rmdefault}{\mddefault}{\updefault}\special{ps: gsave 0 0 0 setrgbcolor}$R_{1,123}(X)$\special{ps: grestore}}}}
\put(8191,-2131){\makebox(0,0)[lb]{\smash{\SetFigFont{6}{7.2}{\rmdefault}{\mddefault}{\updefault}\special{ps: gsave 0 0 0 setrgbcolor}$R_{123}^{1}(X)$\special{ps: grestore}}}}
\put(8101,-3706){\makebox(0,0)[lb]{\smash{\SetFigFont{6}{7.2}{\rmdefault}{\mddefault}{\updefault}\special{ps: gsave 0 0 0 setrgbcolor}$R_{123}^{123}(X)$\special{ps: grestore}}}}
\end{picture}